\documentclass[english,a4paper]{article}
\usepackage[english]{babel}
\usepackage[latin1]{inputenc}
\usepackage{amsmath,amsfonts}
\usepackage{mathrsfs}
\newtheorem{Theorem}{Theorem}[section]

\newtheorem{Remark}[Theorem]{Remark}
\newtheorem{Proposition}[Theorem]{Proposition}
\newtheorem{Corollary}[Theorem]{Corollary}
\numberwithin{equation}{section}

\title{A Proposal of Multigrid Methods for Hermitian Positive Definite Linear Systems enjoying an order relation}
\author{Stefano Serra-Capizzano\thanks{Dipartimento di Scienza e alta Tecnologia,
Universit\`a dell'Insubria - sede di Como, via Valleggio, 11,
22100 Como, Italy. E-mail:stefano.serrac@uninsubria.it} \and
Cristina Tablino-Possio\thanks{Dipartimento di Matematica e
Applicazioni, Universit\`a di Milano Bicocca, Via Cozzi, 53, 20125
Milano, Italy. E-mail:cristina.tablinopossio@unimib.it.} }
\begin{document}
\maketitle
\begin{abstract}
\noindent
Given a multigrid procedure for linear systems with coefficient matrices $A_n$,
we discuss the optimality of a related multigrid procedure with the same smoother
and the same projector, when applied to properly related
algebraic problems with coefficient matrices $B_n$: we assume that
both $A_n$ and $B_n$ are positive definite with $A_n\le \vartheta B_n$, for some positive $\vartheta$
independent of $n$.  In this context we prove the Two-Grid method
optimality. We apply this elementary strategy for designing
a multigrid solution for modifications of multilevel structured (Toeplitz,
circulants, Hartley, sine ($\tau$ class) and cosine algebras)
linear systems, in which the coefficient matrix is banded in a
multilevel sense and Hermitian positive definite.
In such a way, several linear systems arising
from the approximation of integro-differential equations with
various boundary conditions can be efficiently solved in linear
time (with respect to the size of the algebraic problem). Some
numerical experiments are presented and discussed, both with
respect to Two-Grid and multigrid procedures.
\end{abstract}
\section{Introduction and notations}
In this note we first we discuss some theoretical issues
concerning the Two-Grid Method
(TGM) in terms of the algebraic multigrid theory provided
by Ruge and St\"uben \cite{RStub}. More precisely, we prove that the
proposed TGM is optimally convergent with a convergence rate
independent of the dimension for a given sequence of linear systems
$\{B_ny_n=c_n\}_n$ (with size of $B_n$ increasing function of $n$)
with uniformly bounded Hermitian positive
definite matrix sequence $\{B_n\}_n$, under the assumption that such
TGM is optimal for $\{A_nx_n=b_n\}_n$ with a given Hermitian
positive definite matrix sequence $\{A_n\}_n$ related to $\{B_n\}_n$
by means of a simple order relation. More in detail, we require that there exist
$\vartheta,M>0$ independent of $n$ such that
$A_n\le \vartheta B_n$, with $\|B_n\|_2\le M$, for $n\ge \bar{n}$
(for a discussion on comparison theorems see also \cite{falgout}).
When we write `the same TGM' we mean that the same smoother and the
same projection are used.
Concerning notation and definitions, the following
terminology will be used. By $\|\cdot\|_2$ we denote the Euclidean
norm on $\mathbb{C}^m$ and the associated induced matrix norm over
$\mathbb{C}^{m\times m}$. If $X$ is Hermitian positive definite,
then its square root obtained via the Schur decomposition is well
defined and positive definite. As a consequence we can set
$\|\cdot\|_{X}=\|X^{1/2}\cdot\|_{2}$ the Euclidean norm weighted by
$X$ on $\mathbb{C}^m$, and the associated induced matrix norm. In
addition, the notation $X\le Y$, with $X$ and $Y$ Hermitian
matrices, means that $Y-X$ is nonnegative definite. Furthermore, we
write that $\{X_n\}_n$
is a uniformly bounded matrix sequence if there exists
$M>0$ independent of $n$ such that $\|X_n\|_2\le M$, for $n$ large
enough. Finally, for absolute constant we mean a constant independent of the size
parameter $n$, and when we write that a relation depending on $n$ is true for $n$ large enough we mean
that there exists a value $\bar{n}$ such that the relation is valid for any $n\ge \bar{n}$.
A matrix is called $d$-level, $d>1$, if it is a block matrix
where each block is a $(d-1)$-level matrix. For $d=1$ we find a
standard generic matrix. Furthermore a matrix is banded in $d$-level
sense if it is a block banded matrix where each block is a banded
$(d-1)$-level matrix. For $d=1$ we find a standard banded matrix.
\newline
Here we report several cases interesting in applications or of computational
interest where the key assumption  $A_n\le \vartheta B_n$ is fulfilled with
some positive $\vartheta$ independent of $n$. Among them, as a case study,
we may consider matrices $B_n$ with $B_n=A_n+\Theta_n$,
where $A_n$ is structured, positive definite,  ill-conditioned,
and for which an optimal multigrid algorithm is already available,
and where $\Theta_n$ is an indefinite band correction not
necessarily structured; moreover, we require that $\{A_n\}_n$ and
$\{B_n\}_n$ are uniformly bounded and that  $A_n+\Theta_n$ is
still positive definite and larger than $A_n/\vartheta$ for some
$\vartheta>0$ independent of $n$. \newline For instance such a
situation is encountered when
dealing with standard Finite Difference (FD)
approximations of the problem
$$
{\cal L} u  =  -\Delta  u(x) + \mu(x) u(x)
  =  h(x),\ x \in \Omega,
$$
where $\mu(x)$ and $h(x)$ are given bounded functions,
$\Omega=(0,1)^d$, $d\ge 1$, and with Dirichlet, periodic or
reflective boundary conditions (for a discussion on various
boundary conditions see \cite{NCT,model-tau}). For specific
contexts where structured-plus-diagonal problems arise refer to
\cite{mehrmann,ng-diag-ci} and \cite{BGST}, when considering also
a convection term in the above equation. However, the latter is
just an example chosen for the relevance in applications, but the
effective range of applicability of our proposal is indeed much
wider.\newline
The numerical experimentation suggests that an
optimal convergence rate should hold for the MGM as well. Here,
for MGM algorithm, we mean the simplest (and less expensive)
version of the large family of multigrid methods, i.e., the
V-cycle procedure: for a brief description of the TGM and of the
V-cycle algorithms we refer to Section \ref{sec:algmul}, while an
extensive treatment can be found in \cite{hack} and especially in
\cite{Oost}. Indeed, we remark that in all the considered
examples the MGM is optimal in the sense that (see \cite{AxN}):
\begin{itemize}
\item[\textbf{a.}] the observed number of iterations is constant
with respect to the size of the algebraic problem;
\item[\textbf{b.}] the cost per iteration (in terms of arithmetic
operations) is just linear as the size of the algebraic problem.
\end{itemize}
Nevertheless, it is worth stressing that the theoretical extension of the
optimality result to the MGM is still an open
question, while the analysis can be easily extended to a fixed number of levels.
\newline

We should mention that an extensive literature has treated the
numerical solution of structured linear systems of large dimensions
\cite{CN}, by means of preconditioned iterative solvers. However, as
well known in the multilevel setting, the most popular matrix
algebra preconditioners cannot work in general (see \cite{TCS} and
references therein), while the multilevel structures often are the
most relevant in practical applications.
\newline
Therefore, quite recently, more attention has been focused (see
\cite{AD,CSun2,FS2,Hu,Sun,mcircoST-weighted}) on the multigrid
solution of multilevel structured (Toeplitz, circulants, Hartley,
sine ($\tau$ class) and cosine algebras) linear systems, in which
the coefficient matrix is banded in a multilevel sense and
Hermitian positive definite. The reason is due to the fact that
these techniques are very efficient, the total cost for reaching
the solution within a preassigned accuracy being linear as the
dimensions of the involved linear systems.
%
\newline
The paper is organized as follows. In Section \ref{sec:algmul} we
report the standard  TGM and MGM algorithms and we write explicitly the
related iteration matrices.
In Section \ref{sec:extension} we first recall some
classical results related to the algebraic TGM convergence analysis
and then we prove the optimal rate of
convergence of the proposed TGM, under the order relation previously
introduced. The MGM case is briefly discussed at the end of the
section. In Section \ref{sec:order} we report different interesting
situations in which the key order relation is satisfied and, as a
case study, we consider the discrete Laplacian-plus-diagonal systems.
Finally in Section \ref{sec:numexp} we discuss several numerical experiments, while
Section \ref{sec:finrem} deals with further considerations
concerning future work and perspectives.

\section{Two-grid and Multigrid iterations}\label{sec:algmul}
Let $n_0$ be a positive $d$-index, $d\ge 1$, and let $N(\cdot)$ be
an increasing function with respect to $n_0$. In devising a TGM
and a MGM for the linear system
$A_{n_0} x_{n_0}=b_{n_0}$,
where $A_{n_0} \in \mathbb{C}^{N({n_0})\times N({n_0})}$  and
$x_{n_0},b_{n_0} \in \mathbb{C}^{N({n_0})}$, the ingredients below
must be considered.\newline
Let $n_1 < n_0$ (componentwise) and let $p_{n_0}^{n_1}\in
\mathbb{C}^{N({n_0})\times N({n_1})}$ be a given full-rank matrix.
In order to simplify the notation, in the following we will refer to
any multi-index $n_s$ by means of its subscript $s$, so that, e.g.
$A_s:=A_{n_s}$, $b_s:=b_{n_s}$, $p_s^{s+1}:=p_{n_s}^{n_{s+1}}$, etc.
With the same notations, a class of stationary iterative methods of the
form
$x_s^{(j+1)}=V_s x_s^{(j)}+ \tilde b_s$
is also considered in such a way that ${\mathcal
Smooth}(x_s^{(j)},b_s,V_s,\nu_s)$ denotes the application of this
rule $\nu_s$ times, with $\nu_s$ positive integer number, at the
dimension corresponding to the index $s$. \newline
Thus, the solution of the
linear system $A_{n_0} x_{n_0}=b_{n_0}$ is obtained by applying
repeatedly the TGM iteration, where the $j^{\mathrm th}$ iteration\newline
\[
{x}_0^{(j+1)} = {\mathcal TGM}({x}_0^{(j)}, {b}_0,
A_0,V_{0,\mathrm{pre}},\nu_{0,\mathrm{pre}},
V_{0,\mathrm{post}},\nu_{0,\mathrm{post}})
\]
is defined by the following algorithm \cite{hack}:\newline
\[
   \begin{array}{l}
     \multicolumn{1}{c}
{y_0:={\mathcal TGM}({x}_0, {b}_0, A_0,V_{0,\mathrm{pre}},\nu_{0,\mathrm{pre}},V_{0,\mathrm{post}},\nu_{0,\mathrm{post}})} \\
     \hline \\
    {\bf 1.}\ \ \tilde x_0:={\mathcal Smooth}(x_0,b_0,V_{0,\mathrm{pre}},\nu_{0,\mathrm{pre}}) \\
    {\bf 2.}\ \  r_{0}:= b_0- A_0 \tilde x_0\\
    {\bf 3.}\ \  r_{1}:=(p_0^1)^H r_{0}     \\
    {\bf 4.}\ \  \mbox{Solve\ } A_{1}y_1=r_{1}, \textrm{ with }  A_{1}:=(p_0^1)^HA_0 p_0^1 \\
    {\bf 5.}\ \  \tilde y_0:=\tilde x_0+p_0^1 y_1\\
    {\bf 6.}\ \  y_0:={\mathcal Smooth}(\tilde y_0,b_0,V_{0,\mathrm{post}},\nu_{0,\mathrm{post}})
   \end{array}
\]
Steps $1.$ and $6.$ concern the application of
$\nu_{0,\mathrm{pre}}$ steps of the \emph{pre-smoothing} (or
\emph{intermediate}) iteration and of $\nu_{0,\mathrm{post}}$ steps
of the \emph{post-smoothing} iteration, respectively. Moreover,
steps $2.\rightarrow 5.$ define the so called \emph{coarse grid
correction}, that depends on the projection operator $(p_0^1)^H$. In
such a way, the TGM iteration represents a classical stationary
iterative method whose iteration matrix is given by
\begin{equation}\label{explicit:tgm-exact}
TGM_0= V_{0,\mathrm{post}}^{\nu_{0,\mathrm{post}}}
\ 
  CGC_0
\ 
V_{0,\mathrm{pre}}^{\nu_{0,\mathrm{pre}}},
\end{equation}
where $CGC_0=   I_0-p_0^1
            \left[ (p_0^1)^H A_0 p_0^1  \right]^{-1}
   (p_0^1)^H A_0$
denotes the coarse grid correction iteration matrix. The names
\emph{intermediate} and \emph{smoothing iteration} used above refer
to the multi-iterative terminology \cite{Smulti}: we say that a
method is multi-iterative if it is composed by at least two distinct
iterations. The idea is that these basic components should have
complementary spectral behaviors so that the whole procedure is
quickly convergent (for details see \cite{Smulti} and Sections 7.2
and 7.3 in \cite{mcirco}).
Notice
that in our setting of Hermitian positive definite and uniformly
bounded sequences, the subspace where $A_0$ is ill-conditioned
corresponds to the subspace in which $A_0$ has small eigenvalues.
\newline
Starting from the TGM, the MGM is easily introduced. Indeed, the main
difference with respect to the TGM is as follows: instead of
solving directly the linear system with coefficient matrix $A_1$,
we apply recursively the projection strategy so obtaining a
{\em multigrid method}. \newline
Let us use the Galerkin
formulation and let $n_0>n_1>\ldots
>n_l>0$, with $l$ being the maximal number of recursive calls
and with $N(n_s)$ being the corresponding matrix sizes.
\newline
The corresponding multigrid method generates the approximate
solution at the $j^{\mathrm th}$ iteration
\[
{x}_0^{(j+1)} = {\mathcal MGM}(0, {x}_0^{(j)}, {b}_0, A_0,
V_{0,\mathrm{pre}},{\nu}_{0,\mathrm{pre}}, V_{0,\mathrm{post}},
{\nu}_{0,\mathrm{post}})
\]
according to the following algorithm:
\[
    \begin{array}{c}
        {y}_{s}:={\mathcal MGM}(s,{x}_{s},{b}_{s}, A_{s}, V_{{s},\mathrm{pre}},{\nu}_{s,\mathrm{pre}},V_{{s},\mathrm{post}},{\nu}_{s,
        \mathrm{post}}) \\
        \hline                                                        \\[-2.5mm]
        \begin{array}{l@{}l@{}l@{}l}
            \mathrm{If \  }s=l \ & \mathrm{then\ } & \multicolumn{2}{@{}l}{{\mathcal Solve}(A_{s}{y}_{s}={b}_{s})}            \\
            \                       &\mathrm{else } & {\bf 1.}\ \ \ {\tilde{x}}_{s}     &:= {\mathcal Smooth}\left({x}_{s}, {b}_{s}, V_{{s},
            \mathrm{pre}}, \nu_{s,\mathrm{pre}} \right) \\
            \                       &            & {\bf 2.}\ \ \ {r}_{s}      &:= {b}_{s} - A_{s}{\tilde{x}}_{s}  \\
            \                       &            & {\bf 3.}\ \ \ {r}_{s+1}   &:= (p_{s}^{s+1})^H {r}_{s}   \\
            \                       &            & \multicolumn{2}{@{}l} {{\bf 4.} \quad {y}_{s+1}\!\!:=
    {\mathcal MGM}(s+1,{0}_{s+1},{r}_{s+1},A_{s+1},\!V_{{s+1},\mathrm{pre}},{\nu}_{s+1,\mathrm{pre}},}
\\
            \                       &            &
            \multicolumn{2}{@{}l}{ \hskip 3cm   \!V_{{s+1},\mathrm{post}},{\nu}_{s+1,\mathrm{post}})}  \\
            \ &  & {\bf 5.}\ \ \ \tilde{y}_{s} &:= {\tilde{x}}_{s}+  p_{s}^{s+1}{y}_{s+1} \\
            \  &  & {\bf 6.}\ \ \ {{y}}_{s}  &:= {\mathcal Smooth}\left( \tilde{y}_{s}, {b}_{s}, V_{{s},\mathrm{post}}, \nu_{s,
            \mathrm{post}} \right)
         \end{array}
    \end{array}
\]
where the matrix $A_{s+1} := (p_{s}^{s+1})^H A_{s} p_{s}^{s+1}$ is
more profitably computed in the so called \emph{pre-computing
phase}. \newline Since the multigrid is again a linear fixed-point
method, we can express ${x}_0^{(j+1)}$ as $MGM_{0}{x}_0^{(j)}$
$+(I_0-MGM_{0})A_0^{-1}{b}_0$,  where the iteration matrix
$MGM_{0}$ is recursively defined  according to the following rule
(see \cite{Oost}):
\begin{equation}\label{eq:MGM}\left\{\!
  \begin{array}{@{}l@{\;}c@{\;}l}
    MGM_{l} & = & O,      \\[3mm]
    MGM_{s} & = & V_{s,\mathrm{post}}^{\nu_{s,\mathrm{post}}} \left[
      I_{s}\!-\! p_{s}^{s+1} \!\!\left(
      I_{s+1}\!-\!MGM_{s+1}\right) \!\!
      A_{s+1}^{-1}(p_{s}^{s+1})^H\, A_{s}
      \right]  V_{ s,\mathrm{pre}}^{ \nu_{s,\mathrm{pre}}}, \\
    &  & \qquad\qquad \qquad\qquad \qquad\qquad \qquad\qquad
         \qquad  s=0, \dots, l-1,
  \end{array}
\right.\!
\end{equation}
and with $MGM_s$ and $MGM_{s+1}$ denoting the iteration matrices
of the multigrid procedures at two subsequent levels, $s=0, \dots,
l-1$.
At the last recursion level $l$, the linear system is solved by a
direct method and hence it can be interpreted as an iterative
method converging in a single step: this motivates the
initial condition  $MGM_{l} = O$.\newline
By comparing the TGM and MGM, we observe that  the coarse grid
correction operator $CGC_s$ is replaced by an approximation, since
the matrix $A_{s+1}^{-1}$ is approximated by
$\left(I_{s+1}-\!MGM_{s+1}\right)A_{s+1}^{-1}$ as implicitly
described in (\ref{eq:MGM}) for $s=0,\ldots,l-1$. In this way step
$4.$, at the highest level $s=0$, represents an approximation of
the exact solution of step $4.$ displayed in the TGM algorithm
(for the matrix analog compare (\ref{eq:MGM}) and
(\ref{explicit:tgm-exact})).
Finally, for $l=1$  the MGM
reduces to the TGM if ${\mathcal Solve}(A_1{y}_1={b}_1)$ is
${y}_1=A_1^{-1}{b}_1$.
\section{Discussion and extension of known convergence results} \label{sec:extension}
In the algebraic multigrid theory some relevant convergence
results are due to Ruge and St\"uben \cite{RStub}. In fact, they
provide the main theoretical tools to which we refer in order to
prove our subsequent convergence results. More precisely,
by referring to the work of Ruge and St\"uben \cite{RStub}, we
will consider Theorem 5.2 therein  in its original form and in the
case where both pre-smoothing and post-smoothing iterations are
performed. In the following all the constants 
$\alpha_\mathrm{pre}$, $\alpha_\mathrm{post}$, and $\beta$ are
required to be independent of the actual dimension in order to
ensure a TGM convergence rate independent of the size of the
algebraic problem.
\begin{Theorem}\label{teo:TGMconv}{\rm \cite{RStub}}
Let $A_0$ be a Hermitian positive definite matrix of size
$N(n_0)$, let $p_0^{1}\in \mathbb{C}^{N(n_0)\times N(n_1)}$, $n_0
>n_1$, be a given full-rank matrix and let $V_{0,\mathrm{post}}$ be
the post-smoothing iteration matrix.\newline Suppose that there
exists $\alpha_\mathrm{post}>0$, independent of $n_0$, such that
for all $x\in \mathbb{C}^{N(n_0)}$
\begin{eqnarray}
\|V_{0,\mathrm{post}}x\|_{A_0}^2 &\le&
\|x\|_{A_0}^2-\alpha_\mathrm{post}\ \|x\|_{A_0D_0^{-1}A_0}^2,
\label{hp:1post-0}
\end{eqnarray}
where $D_0$ is the diagonal matrix formed by the diagonal entries
of $A_0$. \newline Assume, also, that there exists $\beta>0$,
independent of $n_0$, such that for all $x\in \mathbb{C}^{N(n_0)}$
\begin{equation}\label{hp:2-0}
\min_{y\in \mathbb{C}^{N(n_{1})} } \| x -p_0^1 y \|_{D_0}^2 \le
\beta \ \| x \|_{A_0}^2.
\end{equation}
Then, $\beta \ge \alpha_\mathrm{post}$ and $\ \|TGM_0\|_{A_0} \le
    \sqrt{1-\alpha_\mathrm{post}/\beta}<1$.
\end{Theorem}
\begin{Theorem}\label{teo:TGMconv-pre-post}
Let $A_0$ be a Hermitian positive definite matrix of size
$N(n_0)$, let $p_0^{1}\in \mathbb{C}^{N(n_0)\times N(n_1)}$, $n_0
>n_1$, be a given full-rank matrix and let $V_{0,\mathrm{pre}}$, $V_{0,\mathrm{post}}$ be
the pre-smoothing and post-smoothing iteration matrices,
respectively.\newline Suppose that there exist
$\alpha_\mathrm{pre}$, $\alpha_\mathrm{post}>0$, independent of
$n_0$, such that  for all $x\in \mathbb{C}^{N(n_0)}$
\begin{eqnarray}
\|V_{0,\mathrm{pre}} x\|_{A_0}^2 &\le&
\|x\|_{A_0}^2-\alpha_\mathrm{pre}\
\|V_{0,\mathrm{pre}}x\|_{A_0D_0^{-1}A_0}^2, \label{hp:1pre-1} \\
\|V_{0,\mathrm{post}}x\|_{A_0}^2 &\le&
\|x\|_{A_0}^2-\alpha_\mathrm{post}\ \|x\|_{A_0D_0^{-1}A_0}^2,
\label{hp:1post-1}
\end{eqnarray}
where $D_0$ is the diagonal matrix formed by the diagonal entries
of $A_0$. \newline Assume, also, that there exists $\beta>0$,
independent of $n_0$, such that for all $x\in \mathbb{C}^{N(n_0)}$
\begin{equation}\label{hp:2-1}
\| CGC_0 x \|_{A_0}^2 \le \beta \ \| x \|_{A_0D_0^{-1}A_0}^2.
\end{equation}
Then, $\beta \ge \alpha_\mathrm{post}$ and
\begin{equation}\label{tesi-1}
\|TGM_0\|_{A_0} \le
    \sqrt{ \frac{1-\alpha_\mathrm{post}/\beta}{1+\alpha_\mathrm{pre}/\beta}}<1.
\end{equation}
\end{Theorem}
\begin{Remark}\label{rem:In} Theorems \emph{\ref{teo:TGMconv}} and
\emph{\ref{teo:TGMconv-pre-post}} still hold if the diagonal
matrix $D_0$ is replaced by any Hermitian positive matrix $X_0$
\emph{(}see e.g. {\rm \cite{ADS}}\emph{)}. More precisely, $X_0=I$ could be a
proper choice for its simplicity.
\end{Remark}
\begin{Remark}\label{rem:approx-cond}
For reader convenience, the essential steps of the proof of
Theorems \emph{\ref{teo:TGMconv}} and
\emph{\ref{teo:TGMconv-pre-post}} are reported in Appendix
\ref{sez:appendix}, where relations \emph{(\ref{hp:1post-0})} and
\emph{(\ref{hp:1pre-1})} are called post-smoothing and
pre-smoothing property, respectively, and the relation
\emph{(\ref{hp:2-1})} is called approximation property. In this
respect, we notice that the approximation property deduced by
using \emph{(\ref{hp:2-0})} holds only for vectors belonging to
the range of $CGC_0$; 
conversely the approximation property described in \emph{(\ref{hp:2-1})} is
unconditional, i.e., it is satisfied for all $x\in
\mathbb{C}^{N(n_0)}$.
\end{Remark}
In this paper we are interested in the multigrid solution of
special linear systems of the form
\begin{equation}\label{bn-system}
B_n x=b, \quad B_n\in \mathbb{C}^{N(n)\times N(n)}, \ x,b
\in\mathbb{C}^{N(n)}
\end{equation}
with $\{B_n\}_n$ Hermitian positive definite uniformly bounded
matrix sequence, $n$ being a positive $d$-index, $d\ge 1$ and
$N(\cdot)$ an increasing function with respect to it.  More
precisely, we assume that there exists $\{A_n\}_n$ Hermitian
positive definite matrix sequence such that some order relation is
linking $\{A_n\}_n$ and $\{B_n\}_n$,
for $n$ large enough. We suppose also that an optimal algebraic
multigrid method is available for the solution of the systems
\begin{equation}\label{an-system}
A_n x=b, \quad A_n\in \mathbb{C}^{N(n)\times N(n)},\ x,b
\in\mathbb{C}^{N(n)}.
\end{equation}
We ask wether the algebraic TGM and MGM considered for the systems
(\ref{an-system}) are effective also for the systems
(\ref{bn-system}), i.e., when considering 
the very same projectors. Since it is well-known that a very
crucial role in algebraic MGM is played by the choice of projector operator,
the quoted choice will give rise to a relevant simplification. The
results pertain to the convergence analysis of the TGM and MGM: we
provide a positive answer for the TGM case and we only discuss the
MGM case, which is substantially more involved.
%
%
\subsection{TGM convergence and optimality: theoretical
results}\label{tgm-opt}
In this section we give a theoretical analysis of the TGM  in
terms of the algebraic multigrid theory due to Ruge and St\"uben
\cite{RStub} according to Theorem \ref{teo:TGMconv}.
\begin{Proposition}\label{prop:TGMconv-A}
Let $\{A_n\}_n$ be a matrix sequence with $A_n$ Hermitian positive
definite matrices and  let $p_0^1\in \mathbb{C}^{N(n_0)\times
N(n_1)}$ be a given full-rank matrix for any $n_0>0$ such that
there exists $\beta_A>0$ independent of $n_0$ so that for all
$x\in \mathbb{C}^{N(n_0)}$
\begin{equation}\label{hp:2A}
\min_{y\in\mathbb{C}^{N(n_1)} } \| x -p_0^1 y \|_{2}^2 \le \beta_A
\| x \|_{A_0}^2.
\end{equation}
Let $\{B_n\}_n$ be another matrix sequence, with $B_n$ Hermitian
positive definite matrices, such that $A_n\le \vartheta B_n$, for
$n$ large enough, with $\vartheta>0$ absolute constant. Then, by setting $\beta_B=\beta_A \vartheta$ it also holds that, for all $x\in
\mathbb{C}^{N(n_0)}$ and $n_0$ large enough, we have
\begin{equation}\label{hp:2B}
\min_{y\in \mathbb{C}^{N(n_1)} } \| x -p_0^1 y \|_{2}^2 \le
\beta_B \| x \|_{B_0}^2.
\end{equation}
\end{Proposition}
{\bf Proof.}\ From (\ref{hp:2A}) and from the assumptions on the
order relation, we deduce that for all $x\in \mathbb{C}^{N(n)}$
$\min_{y\in \mathbb{C}^{N(n_1)}} \|x_-p_0^1 y\|_{2}^2  \le
\beta_A \|x\|_{A_0}^2 \le  \vartheta \beta_A \|x\|_{B_0}^2$,
i.e., taking into account Remark \ref{rem:In}, the hypothesis
(\ref{hp:2-0}) of Theorem \ref{teo:TGMconv} is fulfilled  for
$\{B_n\}_n$ too, with constant $\beta_B=\beta_A\vartheta$, by
considering the very same projector $p_0^1$ considered for
$\{A_n\}_n$. \hfill $\bullet$ 5\newline
\newline
Thus, the convergence result in Theorem \ref{teo:TGMconv} holds
true also for the matrix sequence $\{B_n\}_n$, if we are able to
guarantee also the validity of condition (\ref{hp:1post-0}). It is
worth stressing that  in the case of Richardson smoothers such
topic is not related to any partial ordering relation connecting
the Hermitian matrix sequences $\{A_n\}_n$ and $\{B_n\}_n$. In
other words, given a partial ordering between $\{A_n\}_n$ and
$\{B_n\}_n$, inequalities (\ref{hp:1post-0}), (\ref{hp:1pre-1}),
and (\ref{hp:1post-1}) with $\{B_n\}_n$ instead of $\{A_n\}_n$ do
not follow from (\ref{hp:1post-0}), (\ref{hp:1pre-1}), and
(\ref{hp:1post-1}) with $\{A_n\}_n$, but they have to be proved
independently. See Proposition 3 in \cite{AD} for the analogous claim in the case
of $\nu_{\mathrm{pre}}, \nu_{\mathrm{post}} >0$.
\begin{Proposition}\label{prop:TGMconv-S}
Let $\{B_n\}_n$ be an uniformly bounded matrix sequence, with
$B_n$ Hermitian positive definite matrices.
For any $n_0>0$, let $V_{0,
\mathrm{pre}}=I_0-\omega_{\mathrm{pre}} B_0$, $V_{0,
\mathrm{post}}=I_0-\omega_{\mathrm{post}} B_0$ be the
pre-smoothing and post-smoothing iteration matrices, respectively
considered in the $TGM$ algorithm, with $\omega_{\mathrm{pre}},\omega_{\mathrm{post}}\in (0, 2/M)$, $M= \textrm{sup}_{n_0>0} \rho(B_0)$. Then, there exist
$\alpha_{B,\mathrm{pre}}$, $\alpha_{B,\mathrm{post}}>0$
independent of $n_0$ such that for all $x\in \mathbb{C}^{N(n_0)}$
\begin{eqnarray}
\|V_{0,\mathrm{pre}} x\|_{B_0}^2 &\le&
\|x\|_{B_0}^2-\alpha_{B,\mathrm{pre}}
\|V_{0,\mathrm{pre}}x\|_{B_0^2}^2,  \label{hp:1preB} \\
\|V_{0,\mathrm{post}} x\|_{B_0}^2 &\le&
\|x\|_{B_0}^2-\alpha_{B,\mathrm{post}} \|x\|_{B_0^2}^2.
\label{hp:1postB}
\end{eqnarray}
\end{Proposition}
{\bf Proof.}\ Relation (\ref{hp:1postB}) is equivalent to the
existence of an absolute  constant
$\alpha_{B,\mathrm{post}}>0$ such that
$(I_0-\omega_{\mathrm{post}}B_0)^2 B_0\le B_0-
\alpha_{B,\mathrm{post}} B_0^2$,
i.e.,
$\omega_{\mathrm{post}}^2 B_0 -2\omega_{\mathrm{post}} I_0 \le
-\alpha_{B,\mathrm{post}} I_0$. \newline
The latter is equivalent to require that the inequality
$\alpha_{B,\mathrm{post}} \le \omega_{\mathrm{post}}(2-
\omega_{\mathrm{post}}\lambda) $ is satisfied for any eigenvalue
$\lambda$ of the Hermitian matrix $B_{0}$ with
$\alpha_{B,\mathrm{post}}>0$ independent of $n_0$. Now, let $[m,M]$
be the smallest interval containing the topological closure of the union over
all $n$ of all the eigenvalues of $B_n$. Then it is enough to
consider
$
\alpha_{B,\mathrm{post}} \le \omega_{\mathrm{post}} \min_{\lambda
\in [m,M]}(2- \omega_{\mathrm{post}}
\lambda)=\omega_{\mathrm{post}} (2- \omega_{\mathrm{post}} M)$,
where the condition $\omega_{\mathrm{post}} <  2/M$ ensures
$\alpha_{B,\mathrm{post}}>0$.\newline By exploiting an analogous
technique, in the case of relation (\ref{hp:1preB}), it is
enough to consider
\begin{eqnarray*}
\alpha_{B,\mathrm{pre}} &\le& \omega_{\mathrm{pre}} \min_{\lambda
\in [0,M]}
\frac{\omega_{\mathrm{pre}}(2-\omega_{\mathrm{pre}}\lambda)}{(1-\omega_{\mathrm{pre}}\lambda)^2}\\
&=& \left \{
\begin{array}{lcl}
2\omega_{\mathrm{pre}}&& \textrm{if\ } 0 <\omega_{\mathrm{pre}} \le 3/(2M), \\
&&\\
\displaystyle \frac{\omega_{\mathrm{pre}}
(2-\omega_\mathrm{pre}M)} { (1-\omega_{\mathrm{pre}} M)^2}&&
\textrm{if\ } 3/(2M) \le\omega_{\mathrm{pre}}< 2/M,
\end{array}
\right.
\end{eqnarray*}
where we consider the only interesting case $m=0$, since $m>0$ is
related to the case of well-conditioned systems where the smoother alone (or a basic conjugate gradient method) is an
optimal method so that we do not need the use of a multigrid technique.
\hfill $\bullet$ 
\newline
In this way, according to the Ruge and St\"uben algebraic theory,
we have proved the TGM optimality, that is its convergence
rate independent of the size $N(n)$ of the involved algebraic
problem.
\begin{Theorem}\label{tgm-only-post}
Let $\{B_n\}_n$ be an uniformly bounded matrix sequence, with
$B_n$ Hermitian positive definite matrices. Under the same
assumptions of Propositions \emph{\ref{prop:TGMconv-A}} and
\emph{\ref{prop:TGMconv-S}} the TGM with only one step of
post-smoothing converges to the solution of $B_nx=b$ and its
convergence rate is independent of $N(n)$.
\end{Theorem}
{\bf Proof.} By referring to Propositions \ref{prop:TGMconv-A} and
\ref{prop:TGMconv-S} the claim follows according to Theorem
\ref{teo:TGMconv}, where we consider $D_0=I_0$ by virtue of Remark \ref{rem:In}.
 \hfill $\bullet$ \newline
\ \\
Few remarks are useful in order to understand what happens
when
also a pre-smoothing phase is applied.
\begin{description}
\item{A)} The first observation is that the convergence analysis
can be reduced somehow to the case of only post-smoothing. Indeed,
looking at relation (\ref{explicit:tgm-exact}) and recalling that
the spectra of $AB$ and $BA$ are the same for any pair $(A,B)$ of
square matrices (see \cite{bhatia}), it is evident that $TGM_0=
V_{0,\mathrm{post}}^{\nu_{0,\mathrm{post}}} CGC_0
V_{0,\mathrm{pre}}^{\nu_{0,\mathrm{pre}}}$ has the same spectrum,
and hence the same spectral radius $\rho(\cdot)$, as
$V_{0,\mathrm{pre}}^{\nu_{0,\mathrm{pre}}}V_{0,\mathrm{post}}^{\nu_{0,\mathrm{post}}}CGC_0$,
where the latter represents a TGM iteration with only
post-smoothing. Therefore
\begin{eqnarray*}
\rho(TGM_0)& = &
\rho(V_{0,\mathrm{pre}}^{\nu_{0,\mathrm{pre}}}V_{0,\mathrm{post}}^{\nu_{0,\mathrm{post}}}CGC_0)
\le
\|V_{0,\mathrm{pre}}^{\nu_{0,\mathrm{pre}}}V_{0,\mathrm{post}}^{\nu_{0,\mathrm{post}}}CGC_0\|_{A_0}
\\
& \le & \sqrt{1- \frac{ \tilde\alpha_{\mathrm{post}}}{\beta} }
\end{eqnarray*}
where $\tilde\alpha_{\mathrm{post}}$ is the post-smoothing
constant of the cumulative stationary method  described by the
iteration matrix
$V_{0,\mathrm{pre}}^{\nu_{0,\mathrm{pre}}}V_{0,\mathrm{post}}^{\nu_{0,\mathrm{post}}}$.
\item{B)} Setting $\nu_{0,\mathrm{pre}}=\nu_{0,\mathrm{post}}=1$
and with reference to Item A), we easily deduce that
$\tilde\alpha_{\mathrm{post}}\ge \alpha_{\mathrm{post}}$ where the
latter is the post-smoothing constant related to the sole
post-smoothing method with iteration matrix $V_{0,\mathrm{post}}$.
Furthermore, if the two iteration matrices $V_{0,\mathrm{pre}}$
and $V_{0,\mathrm{post}}$ are chosen carefully, i.e., by taking
into account the spectral complementarity principle, then we can
expect that $\tilde\alpha_{\mathrm{post}}$ is much larger than
$\alpha_{\mathrm{post}}$, so that the TGM with both pre-smoothing
and post-smoothing is much faster than that with only
post-smoothing.
\item{C)} Items A) and B) show that the TGM iteration  with both
pre-smoothing and post-smoothing is never worse than the TGM
iteration with only post-smoothing. Therefore Theorem
\ref{tgm-only-post} implies that the TGM  with both post-smoothing
and pre-smoothing is optimal for systems with matrices $B_n$ under
the same assumptions as in Theorem \ref{tgm-only-post}. \newline
\item{D)} At this point the natural question arises: is it
possible to handle directly assumption (\ref{hp:2-1}), instead of
assumption (\ref{hp:2-0})? As observed in Remark
\ref{rem:approx-cond} these two assumptions are tied up and indeed
they represent the approximation property on the range of $CGC_0$
and the approximation property on the whole $\mathbb{C}^{n_0}$,
respectively. However, from a technical viewpoint, they are very
different and in fact we are unable to state a formal analog of
Proposition \ref{prop:TGMconv-A} by using (\ref{hp:2-1}). More
precisely, for concluding that $\| CGC_0 x \|_{A_0}^2 \le \beta_A
\ \| x \|_{A_0^2}^2$ implies $\| CGC_0 x \|_{B_0}^2 \le \beta_B \
\| x \|_{B_0^2}^2$ with $X_0=I$ as in Remark \ref{rem:In} and with
$\vartheta$, $\beta_A$, $\beta_B$ absolute constants, and $A_n\le
\vartheta B_n$, we would need $X\le Y$, $X,Y\ge 0$ implies $X^2
\le \gamma Y^2$ with some $\gamma$ positive and independent of
$n$. The latter with $\gamma=1$ is the operator monotonicity of
the map $z\rightarrow z^2$ which is known to be false in general
\cite{bhatia}. We should acknowledge that there exist important
subclasses of matrices for which $X\le Y$, $X,Y\ge 0$ implies $X^2
\le \gamma Y^2$. However, this matrix theoretic analysis of
intrinsic interest goes a bit far beyond the scope of the present
paper and will be the subject of future investigations.
\item{E)} Remark \ref{rem:approx-cond} furnishes an interesting
degree of freedom that could be exploited. For instance if we
choose $X_0=A_0$, by assuming suitable order relations between
$\{A_n\}_n$ and $\{B_n\}_n$, then proving that $\| CGC_0 x
\|_{A_0}^2 \le \beta_A \ \| x \|_{A_0}^2$ implies $\| CGC_0 x
\|_{B_0}^2 \le \beta_B \ \| x \|_{B_0}^2$ with $\vartheta$,
$\beta_A$, $\beta_B$ absolute constants, becomes easier, but,
conversely, the study of the pre-smoothing and post-smoothing
properties becomes more involved.
\end{description}
\subsection{MGM convergence and optimality: a discussion}\label{mgm-opt}
In this section we briefly discuss the same question as before,
but in connection with the MGM. First of all, we expect that a
more severe assumption between $\{A_n\}_n$ and $\{B_n\}_n$ has to
be fulfilled in order to infer the MGM optimality for $\{B_n\}_n$
starting from the MGM optimality for $\{A_n\}_n$. The reason is
that the TGM is just a special instance of the MGM when setting
$l=1$.\newline In the TGM setting we have assumed a one side
ordering relation: here the most natural step is to consider a two
side ordering relation, that is to assume that there exist
positive constants $\vartheta_1,\vartheta_2$ independent of $n$
such that $\vartheta_1 B_n\le A_n\le \vartheta_2 B_n$, for every
$n$ large enough. The above relationships simply represent the
spectral equivalence condition for sequences of Hermitian positive
definite matrices. In the context of the preconditioned conjugate
gradient method (see \cite{Axelsson}), it is well known that if
$\{P_n\}_n$ is a given sequence of optimal (i.e., spectrally
equivalent) preconditioners for $\{A_n\}_n$, then $\{P_n\}_n$ is
also a sequence of optimal preconditioners for $\{B_n\}_n$ (see
e.g. \cite{TCS}). The latter fact just follows from the
observation that the spectral equivalence is an {\em equivalence}
relation and hence is transitive. \newline In summary, we have
enough heuristic motivations in order to conjecture that the
spectral equivalence is the correct and needed assumption and, in
reality, the numerical experiments reported in Section
\ref{sec:numexp} give a support to the latter statement.
\newline From a theoretical point of view, as done for the TGM, we
start from the Ruge-St\"uben tools \cite{RStub} in the slightly
modified version contained in Theorem 2.3 in \cite{ADS}, that is
taking into account Remark \ref{rem:In} and, for the sake of
simplicity, we assume no pre-smoothing i.e.,
${\nu}_{\mathrm{pre}}=0$. The matrix inequalities coming from the
assumption (2.9) in \cite{ADS} are very intricate since they involve
simultaneously projector operators and smoothers: whence, it is
customary to split it into the \emph{smoothing property} (relation
(2.11a) in \cite{ADS}) and the \emph{approximation property}
(relation (2.11b) in \cite{ADS}). As usual the smoothing property
does not pose any problem. However, we encounter a serious technical
difficulty in the second inequality, i.e., when dealing with the
approximation property. More precisely, we arrive to compare two
Hermitian projectors, depending on the same $p_s^{s+1}$ with the
first involving $A_{n_s}$ and the second involving $B_{n_s}$.
Unfortunately, they can be compared only in very special and too
restricted cases: the needed assumption would not involve ordering,
but only the fact that the columns of $A_{n_s}^{1/2}p_s^{s+1}$ and
those of $B_{n_s}^{1/2}p_s^{s+1}$ span the same space. \newline
However, as already mentioned, the numerical tests tell us that the
latter difficulty is only a technicality and that the right
assumption should involve spectral equivalence. Furthermore, in such
a context of spectral equivalence, a simpler method could be used:
apply $A_n$ as preconditioner for $B_n$ in a PCG method and solve
the linear systems with coefficient matrix $A_n$ by MGM. Of course,
this approach is simpler to implement, but since several linear
systems have to be solved by MGM, often the flop count proves that
applying the MGM directly is more efficient than using it as solver
for the preconditioner. As already stressed, the problem at the
moment is the lack of theoretical results. Therefore, in future
investigations, other directions and proof techniques have to be
explored.
\section{When the order relation is fulfilled} \label{sec:order}

In the present section we consider examples in which the key order
relation is satisfied: {\bf a)} variable coefficient Laplacians
with Dirichlet or Periodic-Dirichlet BCs vs constant coefficient
Laplacians with the same boundary conditions, {\bf b)} multilevel
Toeplitz matrices with nonnegative symbol having a unique zero at
zero vs sine transform matrices with the same symbol,  {\bf c)}
constant coefficient Laplacians with Dirichlet or
Periodic-Dirichlet BCs vs the same Laplacians plus diagonal
matrices.

\subsection{A case study: discrete Laplacian-plus-diagonal systems} \label{sec:laplacian}
We now consider a specific application of the
results in Section \ref{sec:extension}. More precisely, we apply
a multigrid strategy for solving Laplacian-plus-diagonal linear
systems arising from standard Finite Differences (FD)
discretizations of the problem
\begin{equation} \label{eq:ellp}
\begin{array}{llrlll} {\cal L} u
 & = & -\Delta  u(x) + \mu(x) u(x)
 & = & h(x),
 & x \in \Omega, \\
\end{array}
\end{equation}
where $\mu(x)$ and $h(x)$ are given bounded functions,
$\Omega=(0,1)^d$, $d\ge 1$, and with Dirichlet, periodic or
reflective boundary conditions. Thus, we are facing with a matrix
sequence
\begin{equation} \label{eq:bn}
\{B_n\}_n=\{A_n+D_n\}_n,
\end{equation}
where the structure of the matrix sequence $\{A_n\}_n$ is related
both to the FD discretization and to the type of the boundary
conditions and where $\{D_n\}_n$ is a sequence of
uniformly bounded diagonal matrices, due to the
hypothesis that $\mu(x)$ is bounded. %
\newline
Since a fast TGM and MGM working for the Toeplitz (or $\tau$ - the
$\tau$ class is the algebra associated to the most known sine
transform \cite{BC}) part is well-known \cite{FS1,FS2,CCS,Sun}, we
are in position to apply the tools in the preceding section in order
to show that the same technique works, and with a cost linear as the
dimension, in the context (\ref{eq:bn}) too. In the same way, the
extension of suitable MGM procedures proposed in the case of the
circulant \cite{mcirco}, DCT-III cosine
\cite{mcoseni,mcosine-vcycle}, or $\tau$ \cite{FS1,FS2} algebra, can
be considered according to the corresponding boundary conditions.
Clearly, this case study is just an example relevant in
applications, while the results in Section \ref{sec:extension} are
of much wider generality.\newline
Once more, we want to remark that, unfortunately, there is a gap in
the theory with regard to the MGM, even if the numerical tests
reported in Section \ref{sec:numexp} suggest that the MGM applied to
matrices in $\{B_n=A_n+D_n\}_n$ is optimal under the assumptions
that the same MGM is optimal for $\{A_n\}_n$, $A_n$ symmetric
positive definite matrix, and $\{D_n\}_n$  uniformly bounded matrix
sequence, with $A_n\le \vartheta B_n$ for $n$ large enough and for
some fixed $\vartheta>0$ independent of $n$. Clearly if the
matrices $D_n$ are also nonnegative definite then the constant
$\vartheta$ can be set to $1$. This result can be plainly extended
to the case in which $D_n$ is a (multilevel) banded correction.
\subsection{One-Dimensional case}
According to the FD approximation of (\ref{eq:ellp}) with Dirichlet
boundary conditions, we obtain the matrix sequence
$\{B_n\}_n=\{ A_n+D_n\}_n$
where
$\{A_n\}_n=\{\mathrm{tridiag}_n \ [-1, 2, -1]\}_n$ and
$\{D_n\}_n$ is a sequence of diagonal matrices whose diagonal
 entries $d_{i}^{(n)}$, $i=1,\ldots,n$, are uniformly bounded in modulus by
a constant $M$ independent of $n$. Since
\[
\lambda_{\min}(A_n)=4\sin^2\left({\pi\over
2(n+1)}\right)={\pi^2\over n^2}+O(n^{-3}),
\]
we impose the condition
$n^2\min_{1\le i\le n}d_{i}^{(n)}+\pi^2\ge c$
for some $c>0$ independent of $n$ (we consider only the case
$\min_{1\le i\le n}d_{i}^{(n)}< 0$, since the other is trivial).
Thus, also $\{B_n\}_n$ is an uniformly bounded positive definite
matrix sequence and
\[
B_n\ge A_n+{c-\pi^2\over n^2} I \ge {c\over \pi^2} A_n
\]
so satisfying the crucial assumption $A_n\le \vartheta B_n$ in
Proposition \ref{prop:TGMconv-A} with $\vartheta={\pi^2}/{c}$.
\newline
Let us consider $B_0\in \mathbb{R}^{n_0\times n_0}$, with $1$-index
$n_0>0$. Following \cite{FS1,Sun}, we denote by $T_0^1\in
\mathbb{R}^{n_0\times n_1}$, $n_0=2n_1+1$, the operator such that
\begin{equation}\label{def:tnk-tau}
  (T_0^1)_{i,j}=\left\{ \begin{array}{ll}
            1 & \ \ \mbox{\textrm for }\  i=2j, \ \ j=1,\ldots,n_1,  \\
            0 & \ \ \mbox{\textrm otherwise},
          \end{array} \right.
\end{equation}
and we define a projector $(p_0^1)^H$, $p_0^1\in
\mathbb{R}^{n_0\times n_1}$ as
\begin{equation} \label{def:pn-tau}
p_0^1=\frac{1}{\sqrt 2} P_0 T_0^1, \quad P_0=\ \mathrm{tridiag}_0
\ [1, 2, 1].
\end{equation}
Thus, the basic step in order to prove the TGM optimality result
is reported in the proposition below. It is worth stressing that
the claim refers to a tridiagonal matrix correction, since, under
the quoted assumption, each diagonal correction is projected at
the first coarse level into a tridiagonal correction, while the
tridiagonal structure is kept unaltered in all the subsequent
levels.
\begin{Proposition}\label{prop:lowlev1}
Let $B_0=\ \mathrm{tridiag}_0 \ [-1, 2, -1]+T_0 \in
\mathbb{R}^{n_0\times n_0}$, with $n_0>0$ and $T_0$ being a
symmetric uniformly bounded tridiagonal matrix such that $A_0\le
\vartheta B_0$, with $A_0=\mathrm{tridiag}_0 \ [-1, 2, -1]$ and some $\vartheta>0$.\newline
Let $p_0^1=({1}/{\sqrt 2}) \mathrm{tridiag}_0 \ [1, 2, 1] T_0^1$, with
$n_0=2n_1+1$. Then,
$(p_0^1)^H B_0 p_0^1= \mathrm{tridiag}_1 \ [-1, 2, -1]+T_1$
where $T_1\in \mathbb{R}^{n_1\times n_1}$ is a symmetric uniformly
bounded tridiagonal matrix with $A_1\le \vartheta B_1$,
$A_1=\mathrm{tridiag}_1 \ [-1, 2, -1]$, $B_1=(p_0^1)^H B_0 p_0^1$.
\end{Proposition}
{\bf Proof.}\ For the Toeplitz part refer to \cite{Sun}. For the
tridiagonal part we just need a simple check. In fact, the product
$P_0 T_0 P_0$ is a $7$-diagonal matrix ($P_0$ and $T_0$ are
tridiagonal) and the action of $T_0^1$, on the left and on the
right, selects  even rows and columns so that the resulting
matrix is still tridiagonal. Since $A_0\le \vartheta B_0$,
$A_0=\mathrm{tridiag}_0 \ [-1, 2, -1]$, $B_1=(p_0^1)^H B_0 p_0^1$,
and $A_1=\mathrm{tridiag}_1 \ [-1, 2, -1]=(p_0^1)^H A_0 p_0^1$ it
is evident that $A_1 \le \vartheta B_1$. Finally, the uniform
boundedness is guaranteed by the uniform boundedness of all the
involved matrices. \ \hfill $\bullet$
\begin{Corollary}\label{ver:hp2:1d}
Let $B_0=\ \mathrm{tridiag}_0 \ [-1, 2, -1]+T_0 \in
\mathbb{R}^{n_0\times n_0}$, with $n_0>0$ and $T_0$
symmetric tridiagonal matrix such that $A_0\le
\vartheta B_0$, with $A_0=\mathrm{tridiag}_0 \ [-1, 2, -1]$ and some $\vartheta>0$ independent of $n_0$.
Let $p_0^1 =({1}/{\sqrt 2}) \mathrm{tridiag}_0 \ [1, 2, 1] T_0^1$,
with $n_0=2n_1+1$. Then, there exists $\beta_B>0$ independent of
$n_0$ so that inequality \emph{(\ref{hp:2-0})} holds true.
\end{Corollary}
{\bf Proof.}\ Let $A_0=\ \mathrm{tridiag}_0 \ [-1, 2, -1]\in
\mathbb{R}^{n_0\times n_0}$. Then relation (\ref{hp:2-0}) is
fulfilled with the operator $p_0^1$ defined in (\ref{def:pn-tau})
and with a certain $\beta_A$ independent of $n_0$, as proved in
\cite{Sun}. Moreover, from the assumption we have $A_0\le
\vartheta B_0$ so that Proposition \ref{prop:TGMconv-A} implies
that (\ref{hp:2-0}) holds true for $B_0$ with a constant
$\beta_B=\vartheta \beta_A$ independent of $n_0$. \hfill $\bullet$
\begin{Corollary}\label{ver:hp1:1d}
Let $\{B_n\}_n$ be the sequence such that $B_n=\ \mathrm{tridiag}_n
\ [-1, 2, -1]+T_n \in \mathbb{R}^{n\times n}$ with $T_n$ symmetric
uniformly bounded tridiagonal matrices.
Let $V_{n,\mathrm{pre}}=I_n-\omega_{\mathrm{pre}} B_n$, $V_{n,
\mathrm{post}}=I_n-\omega_{\mathrm{post}} B_n$ be the
pre-smoothing and post-smoothing iteration matrices, with
$\omega_{\mathrm{pre}},\omega_{\mathrm{post}}\in (0, 2/\rho(B_n))$.
Then, there exist $\alpha_{B,\mathrm{pre}}$,
$\alpha_{B,\mathrm{post}}>0$ independent of $n$, so that
inequalities \emph{(\ref{hp:1preB})} and \emph{(\ref{hp:1postB})}
hold true.
\end{Corollary}
{\bf Proof.}\ It is evident that $\{B_n\}_n$ is a sequence of
symmetric positive definite matrices uniformly bounded by $4+M$,
with $\|T_n\|_2\le M$ independent of $n$, so that the thesis
follows by the direct application of Proposition
\ref{prop:TGMconv-S}. \hfill $\bullet$
%
\subsection{Two-Dimensional case}
Hereafter, we want to consider the extension of our TGM and MGM to the
case $d>1$. Due to the discretization process, it is natural, and
easier, to work with $d-$indices $n=(n^{(1)},\ldots , n^{(d)})$,
with $n^{(r)}$ integer positive number, $r=1,\ldots,d$. In this
case the matrix dimension is
$N(n_0)=\prod_{r=1}^d n_0^{(r)}$ and when considering the
projected matrices of size $N(n_1)$ we have that $n_1$ is again a
$d-$index and we assume not only $N(n_1)<N(n_0)$, but also
$n_1<n_0$ componentwise.\newline We discuss in detail the
two-level case, since the $d-$level one is a simple generalization.
Thus, in the two-level case, we are dealing with the matrix
sequence $\{B_n\}_n=\{A_n+D_n\}_n$ where $ \{A_n\}_n=\{\mathrm{tridiag}_{n^{(1)}} \ [-1, 2,
-1]\otimes I_{n^{(2)}} + I_{n^{(1)}}\otimes \
\mathrm{tridiag}_{n^{(2)}} \ [-1, 2, -1]\}_n$
and $\{D_n\}_n$ is a sequence of diagonal matrices whose diagonal
 entries $d_{i}^{(n)}$, $i=1,\ldots,N(n)$, are uniformly bounded in modulus by
a constant $M$ independent of $n$. Since
\[
\lambda_{\min}(A_n) \! = \!4\sin^2\!\! \left({\pi\over
2(n^{(1)}+1)}\right)\!+\!4\sin^2\!\! \left({\pi\over
2(n^{(2)}+1)}\right)\!
 = \! {\pi^2\over [n^{(1)}]^2}+\!{\pi^2\over
[n^{(1)}]^2}+O(\psi^{-3}),
\]
we impose the condition $\psi^2\min_{1\le i\le N(n)}d_{i}^{(n)}+\pi^2\ge c$
for some $c>0$ independent of $n$, with $\psi=\min_j\{n^{(j)}\}$.
Thus, also $\{B_n\}_n$ is an uniformly bounded positive definite
matrix sequence and
\[
B_n\ge A_n+{c-\pi^2\over \psi^2} I \ge {c\over \pi^2} A_n
\]
so satisfying the crucial assumption $A_n\le \vartheta B_n$ in
Proposition \ref{prop:TGMconv-A} with $\vartheta={\pi^2}/{c}$.
\newline
The projector definition can be handled in a natural manner by
using tensorial arguments: $(p_0^1)^H \in \mathbb{R}^{N(n_1)\times N(n_0)}$ is constructed in such a
way that
\begin{eqnarray*}
p_0^1&=&P_0 U_0^1\\
 P_0&=&\ \mathrm{tridiag}_{n_0^{(1)}} \ [1, 2,
1]\otimes \
\mathrm{tridiag}_{n_0^{(2)}} \ [1, 2, 1],\\
U_0^1 &=& T_{0}^{1}(n_0^{(1)})\otimes T_{0}^{1} (n_0^{(2)})
\end{eqnarray*}
with $n_0^{(r)}=2 n_1^{(r)}+1$ and where $T_0^1(n_0^{(r)})\in
\mathbb{R}^{n_0^{(r)}\times n_1^{(r)}}$ is the unilevel matrix
given in (\ref{def:tnk-tau}).
Notice that this is the most trivial extension of the unilevel
projector to the two-level setting and  such a choice is also the
less expensive from a computational point of view: in fact,
$p_0^1=\tau_0((2+2\cos(t_1)(2+2\cos(t_2)))U_0^1$ equals
$[\tau_{n_0^{(1)}}(p(2+2\cos(t_1)))T_0^1(n_0^{(1)})] \otimes
[\tau_{n_0^{(2)}}(p(2+2\cos(t_2)))T_0^1(n_0^{(2)})]$.\newline
The proposition below refers to a two-level tridiagonal correction
for the same reasons as in the unilevel case.
\begin{Proposition}\label{prop:lowlev1bis-2}
Let
$B_0=\mathrm{tridiag}_{n_0^{(1)}} \ [-1, 2, -1]\otimes
I_{n_0^{(2)}} + I_{n_0^{(1)}}\otimes \ \mathrm{tridiag}_{n_0^{(
2)}} \ [-1,\\ 2, -1] +T_0 \in \mathbb{R}^{N(n_0)\times N(n_0)}$,
with $n_0>0$ and $T_0$ being a symmetric uniformly bounded tridiagonal block
matrix with tridiagonal blocks such that $A_0\le \vartheta B_0$, with
$A_0=\mathrm{tridiag}_{n_0^{(1)}} \ [-1, 2, -1]\otimes I_{n_0^{(2)}} +
I_{n_0^{(1)}}\otimes \ \mathrm{tridiag}_{n_0^{(2)}} \ [-1, 2, -1]$ and
some $\vartheta>0$ independent of $n_0$. \newline
Let $p_0^1= (\mathrm{tridiag}_{n_0^{(1)}} \ [1, 2, 1]\otimes \
\mathrm{tridiag}_{n_0^{(2)}} \ [1, 2, 1])
(T_{0}^{1}(n_0^{(1)})\otimes T_{0}^{1}(n_0^{(2)}))$,
with $n_0^{(r)}=2n_1^{(r)}+1$, $r=1,2$. Then, $B_1:=(p_0^1)^H B_0
p_0^1$ coincides with $A_1+T_1$ where $T_1\in
\mathbb{R}^{N(n_1)\times N(n_1)}$ is a symmetric uniformly bounded
tridiagonal block matrix with tridiagonal blocks and where $A_1$ is
a two-level $\tau$ tridiagonal block matrix with tridiagonal blocks
asymptotic to
$\mathrm{tridiag}_{n_1^{(1)}} \ [-1, 2, -1]\otimes I_{n_1^{(2)}} +
I_{n_1^{(1)}}\otimes \ \mathrm{tridiag}_{n_1^{(1)}} \ [-1, 2, -1]$,
so that $A_1\le \vartheta B_1$.
\end{Proposition}
{\bf Proof.}\ For the $\tau$ part refer to \cite{Sun}. For the
two-level banded part just a simple check is required. In fact, the product
$P_0 T_0 P_0$ is a $7$-diagonal block matrix with $7$-diagonal
blocks ($P_0$ is a tridiagonal block matrix with tridiagonal blocks)
and the action of $U_0^1$, on the left and on the right, selects
even rows and columns in even blocks with respect to the rows and
columns, so that the resulting matrix has a tridiagonal block
pattern with tridiagonal blocks. The order relation follows as a
direct consequence  of the assumption $A_0\le \vartheta B_0$ and the
uniform boundedness is implied by the uniform boundedness of all the
involved matrices. \hfill $\bullet$ 
\begin{Corollary}\label{ver:hp2:2d-2}
Let $B_0\!=\! A_0+T_0\! \in \! \mathbb{R}^{\! N(n_0)\times N(n_0)}$ with
$n_0>0$, $A_0 \!=\! \mathrm{tridiag}_{n_0^{(1)}}[-1,\\ 2,
-1]\otimes I_{n_0^{(2)}} + I_{n_0^{(1)}}\otimes \
\mathrm{tridiag}_{n_0^{(2)}} \ [-1, 2, -1]$,
and $T_0$ symmetric tridiagonal block matrix with tridiagonal blocks
such that $A_0\le \vartheta B_0$ for some $\vartheta>0$ independent of $n_0$. Let
$p_0^1= (\mathrm{tridiag}_{n_0^{(1)}} \ [1, 2, 1]\otimes \
\mathrm{tridiag}_{n_0^{(2)}} \ [1, 2, 1]) (T_{0}^{1}(n_0^{(1)})
\otimes T_{0}^{1}(n_0^{(2)}))$,
with $n_0^{(r)}=2n_1^{(r)}+1$, $r=1,2$. Then, there exists
$\beta_B>0$ independent of $n_0$  so that inequality
\emph{(\ref{hp:2-0})} holds true.
\end{Corollary}
{\bf Proof.}\ The proof can be done following the same steps as in
Corollary \ref{ver:hp2:1d}. \hfill $\bullet$ 
\begin{Corollary}\label{ver:hp1:2d-2}
Let $\{B_n \}_n$ be the sequence such that $B_n=A_n+T_n$ with
$A_n = \mathrm{tridiag}_{n^{(1)}} \ [-1, 2, -1]\otimes I_{n^{(2)}}
+ I_{n^{(1)}}\otimes \ \mathrm{tridiag}_{n^{(2)}} \ [-1, 2, -1]$,
and with $T_n$ symmetric uniformly bounded tridiagonal block
matrix with tridiagonal blocks.
%
Let $V_{n,\mathrm{pre}}=I_n-\omega_{\mathrm{pre}} B_n$, $V_{n,
\mathrm{post}}=I_n-\omega_{\mathrm{post}} B_n$ be the
pre-smoothing and post-smoothing iteration matrices, with
$\omega_{\mathrm{pre}},\omega_{\mathrm{post}}\in (0, 2/\rho(B_n))$.
Then, there exist  $\alpha_{B,\mathrm{pre}}$,
$\alpha_{B,\mathrm{post}}>0$ independent of $n$, so that
inequalities \emph{(\ref{hp:1preB})} and \emph{(\ref{hp:1postB})}
hold true.
\end{Corollary}
{\bf Proof.}\ The proof can be worked out as in Corollary
\ref{ver:hp1:1d} since the sequence $\{B_n\}_n$ is uniformly
bounded by $8+M$, with $\|T_n\|_2\le M$ independent of $n$ by
assumption. \ \hfill $\bullet$ 
\section{Numerical Evidences}\label{sec:numexp}
We test our TGM and MGM (standard V-cycle according to Section
\ref{sec:algmul}) for several examples of matrix corrections
$\{D_n\}_n$, $D_n \in \mathbb{C}^{N(n)\times N(n)}$,
$N(n)=\prod_{r=1}^d n^{(r)}$, $d=1,2$.\newline
We will consider nonnegative definite band corrections and
indefinite band corrections. By referring to Section
\ref{sec:laplacian}, the case of nonnegative definite corrections
implies trivially that $A_n\le B_n$ so that the desired constant
is $\vartheta=1$. However, as observed in real-world applications
(see \cite{mehrmann}), the most challenging situation is the one
of indefinite corrections. \newline Concerning nonnegative
definite corrections, the reference set is defined according to
the following notation, in the unilevel and in the two-level
setting, respectively:\newline
\begin{center}
\begin{tabular}{lc ccccc}
& $d_{s}^{(n)}$ & d0& d1& d2& d3 & d4\\ \hline
1D & $s=1,\ldots,N(n)$ & $0$
                     & $\frac{s}{s+1}$
                     & $|\sin(s)|$
                     & $|\sin(s)|\frac{s^2-1}{s^2+1}$
                     & $\frac{s}{N(n)}$ \\
\hline
   & $s=(i-1)n_1+j$  & $0$
                     & $\ \frac{i}{i+1}$
                     & $\ |\sin(i)|$
                     & $\ |\sin(i)|\frac{i^2-1}{i^2+1}$
                     & $\frac{s}{N(n)}$ \\
2D & {\footnotesize $i=1,\ldots,n_1$, $j=1,\ldots,n_2$}   &
                     & $+\frac{j}{j+1}$
                     & $+|\sin(j)|$
                     & $+|\sin(j)|\frac{j^2-1}{j^2+1}$
                     &  \\
\hline
\end{tabular}
\end{center}
\ \\ \ \newline
The case of indefinite corrections is considered in connection with
Laplacian systems with Dirichlet boundary conditions: in that
setting the diagonal entries $d_{s}^{(n)}$ of $D_n$ are generated
randomly. Finally, at the end of the section, higher order differential operators and linear
systems arising from integral equations in image restoration are
considered.
\newline
The aim is to give numerical evidences of the theoretical optimality
results of TGM convergence and to their extension in the case
of the MGM application. \newline
The projectors are properly chosen according to the nature of
structured part, while we will use, in general, the Richardson
smoothing/intermediate iteration step twice in each iteration,
before and after the coarse grid correction, with different values
of the parameter $\omega$.\newline According to the definition,
when considering the TGM, the exact solution of the
system is obtained by using a direct solver 
in the immediately subsequent coarse grid dimension, while, when
considering the MGM, the exact solution of the system is computed by
the same direct solver, when the coarse grid dimension equals $16^d$
(where $d=1$ for the unilevel case and $d=2$ for the two-level
case). 
In all tables we report the numbers of iterations
required for the TGM or MGM convergence, assumed to be reached when
the Euclidean norm of the relative residual becomes less than
$10^{-7}$.
\newline Finally, we stress that the matrices
$A_n$ at every level (except for the coarsest) are never formed
since we need only to store the nonzero Fourier coefficients of the
generating function at every level for matrix-vector
multiplications. For the connection between the entries of a structured matrix ($\tau$, circulant, Toeplitz, etc.) and the
generating function in the context of multigrid techniques refer to
\cite{AD,ADS,Sun,mcirco}.\newline
Thus, besides the $O(N(n))$ operations complexity
of the proposed MGM both with respect to the structured part and
clearly with respect to the non-structured one, the memory
requirements of the structured part are also very low since there
are only $O(1)$ nonzero Fourier coefficients of the generating
function at every level. On the other hand, the projections of the
initial diagonal correction are stored at each level according to
standard sparse matrix techniques during the pre-computing phase.

Concerning the smoothing steps we are aware that Gauss-Seidel (sometimes combined
with one step of conjugate gradient) is often a better choice. Here in order to
check the validity of the approach we consider simpler smoothers like damped Jacobi
or damped Richardson.
\subsection{Discrete Laplacian-plus-diagonal systems}
The numerical tests below refer to convergence results in the case
of matrix sequences arising from the Laplacian discretization, in
the unilevel and in the two-level settings, respectively.
\subsubsection{Dirichlet boundary conditions}
Firstly, we consider the case of Dirichlet boundary conditions so
that the obtained matrix sequence is the Toeplitz/$\tau$ matrix
sequence $\{\tau_n(f)\}_n$ generated by the function
$f(t)=2-2\cos(t)$, $t\in (0,2\pi]$. The projector is defined as in
(\ref{def:tnk-tau}) and (\ref{def:pn-tau}), while the parameters
$\omega$ for the smoothing/intermediate iterations are chosen as
$\omega_{\mathrm{pre}}={2}/{(\|f\|_\infty+\|D_n\|_\infty)}$ and
$\omega_{\mathrm{post}}={1}/{(\|f\|_\infty+\|D_n\|_\infty})$,
with $\nu_\mathrm{pre}=\nu_\mathrm{post}=1$.\newline The results
in the top of Table \ref{table:1d-2d-multigrid} confirm the optimality of the
proposed TGM in the sense that the number of iterations is
uniformly bounded by a constant not depending on the size $N(n)$
indicated in the first column. Moreover, it seems that this claim
can be extended to the MGM convergence. Notice, also, that the
number of iterations is frequently the best possible since it
equals the number of TGM iterations. \newline The case of the
diagonal correction $d4$ deserves special attention: as shown in
the first column, just one pre-smoothing/intermediate and
post-smoothing iteration at each coarse grid level are not
sufficient to ensure the optimality. Moreover, it is enough to
consider a trick that keeps unaltered the $O(N(n))$ computational
cost as proved in \cite{mcirco} (only the multiplicative constant
hidden in the big $O$ can increase): at each projection on a
coarser grid the number of smoothing iterations performed at that
level is increased by a fixed constant $\rho$, i.e., according to
the MGM notation of Section \ref{sec:algmul}, we set
\begin{equation} \label{eq:def-rho}
\nu_{s+1}=\nu_{s}+\rho,\ s=0,\ldots,l-1, \quad \nu_0=1.
\end{equation}
The optimality result in the second column relative to MGM in the
$d4$ case is obtained just by considering $\rho=1$. This
phenomenon is probably due to some inefficiency in considering the
approximation $\| D_n\|_\infty$ in the tuning of the parameter
$\omega_\mathrm{pre}$ and $\omega_\mathrm{post}$. In fact, it is
enough to substitute, for instance, the post-smoother with the
Gauss-Seidel method in order to preserve the optimality also for
$\rho=0$.
\newline Other examples of Toeplitz/$\tau$ linear systems plus
diagonal correction can be found in \cite{nst}, corresponding to
Sinc-Galerkin discretization of differential problems according to
\cite{Lundbook}. Moreover, by using tensor arguments, our results plainly extend to the
two-level setting and the comments concerning the bottom of Table
\ref{table:1d-2d-multigrid} are substantially equivalent as in the unilevel case.
\begin{table}[t!]
\footnotesize
\centering
\begin{tabular}[c]{ccc} 
\multicolumn{3}{c}{{$B_n=\mathrm{tridiag}_{n}[-1, 2,
-1]+$
Diagonal}}\\ 
\begin{tabular}[c]{|c|c|c|c|c|c|}\hline
\multicolumn{6}{|c|}{TGM} \\ \hline
\multicolumn{1}{|c|}{$N(n)$} & \multicolumn{1}{|c|}{d0}
                          & \multicolumn{1}{|c|}{d1}
                          & \multicolumn{1}{|c|}{d2}
                          & \multicolumn{1}{|c|}{d3}
                          & \multicolumn{1}{|c|}{d4}
\\
\multicolumn{1}{|c|}{} & \multicolumn{1}{|c|}{}
                          & \multicolumn{1}{|c|}{}
                          & \multicolumn{1}{|c|}{}
                          & \multicolumn{1}{|c|}{}
                          & \multicolumn{1}{|c|}{}
                          \\
\hline
31   & 2 & 7 & 7 & 7 & 7 \\
63   & 2 & 7 & 8 & 8 & 7 \\
127  & 2 & 7 & 8 & 8 & 7 \\
255  & 2 & 7 & 8 & 8 & 7 \\
511  & 2 & 6 & 8 & 8 & 7 \\
\hline
\end{tabular}
&  &
\begin{tabular}[c]{| c|c|c|c|c|c|c|}\hline
\multicolumn{7}{|c|}{MGM} \\ \hline
\multicolumn{1}{|c|}{$N(n)$} & \multicolumn{1}{|c|}{d0}
                          & \multicolumn{1}{|c|}{d1}
                          & \multicolumn{1}{|c|}{d2}
                          & \multicolumn{1}{|c|}{d3}
                          & \multicolumn{2}{|c|}{d4}
\\
\multicolumn{1}{|c|}{} & \multicolumn{1}{|c|}{}
                          & \multicolumn{1}{|c|}{}
                          & \multicolumn{1}{|c|}{}
                          & \multicolumn{1}{|c|}{}
                          & \multicolumn{1}{|c|}{$\rho\!=\!0$} &  \multicolumn{1}{|c|}{$\rho\!=\!1$}
                          \\
\hline
31 & 2  & 7 & 8 & 8 & 7  &7\\
63  & 7  & 7 & 7 & 7 & 7  &7\\
127 & 8  & 7 & 8 & 8 & 8  &7\\
255 & 8  & 7 & 8 & 8 & 9  &7\\
511 & 8  & 7 & 8 & 8 & 16 &7\\
\hline
\end{tabular}
\end{tabular}
\begin{tabular}[c]{ccc} 
\multicolumn{3}{c}{}\\  
\multicolumn{3}{c}{{$B_n=\mathrm{tridiag}_{n^{(1)}}[-1,
2, -1]\otimes I_{n^{(2)}}+I_{n^{(1)}} \otimes
\mathrm{tridiag}_{n^{(2)}}[-1, 2, -1]
+$Diagonal}}\\  
\begin{tabular}[c]{|c|c|c|c|c|c| }\hline
\multicolumn{6}{|c|}{TGM} \\ \hline
\multicolumn{1}{|c|}{$N(n)$} & \multicolumn{1}{|c|}{d0}
                          & \multicolumn{1}{|c|}{d1}
                          & \multicolumn{1}{|c|}{d2}
                          & \multicolumn{1}{|c|}{d3}
                          & \multicolumn{1}{|c|}{d4}
\\
\multicolumn{1}{|c|}{} & \multicolumn{1}{|c|}{}
                          & \multicolumn{1}{|c|}{}
                          & \multicolumn{1}{|c|}{}
                          & \multicolumn{1}{|c|}{}
                          & \multicolumn{1}{|c|}{}
                          \\
\hline
31$^2$   & 16 & 10 & 13 & 13 & 16 \\
63$^2$   & 16 & 10 & 13 & 13 & 16 \\
127$^2$  & 16 & 10 & 13 & 13 & 16 \\
255$^2$  & 16 & 10 & 13 & 13 & 16 \\
511$^2$  & 16 & 10 & 13 & 13 & 16 \\
\hline
\end{tabular}
& \, &
\begin{tabular}[c]{|c|c|c|c|c|c|c|}\hline
\multicolumn{7}{|c|}{MGM} \\ \hline
\multicolumn{1}{|c|}{$N(n)$} & \multicolumn{1}{|c|}{d0}
                          & \multicolumn{1}{|c|}{d1}
                          & \multicolumn{1}{|c|}{d2}
                          & \multicolumn{1}{|c|}{d3}
                          & \multicolumn{2}{|c|}{d4}
\\
\multicolumn{1}{|c|}{} & \multicolumn{1}{|c|}{}
                          & \multicolumn{1}{|c|}{}
                          & \multicolumn{1}{|c|}{}
                          & \multicolumn{1}{|c|}{}
                          & \multicolumn{1}{|c|}{$\rho\!=\!0$} &  \multicolumn{1}{|c|}{$\rho\!=\!1$}
                          \\
\hline
31$^2$  &16 & 10 & 13 & 13 & 16&16\\
63$^2$  &16 & 10 & 13 & 13 & 17&16\\
127$^2$ &16 & 10 & 12 & 12 & 18&16\\
255$^2$ &16 & 10 & 12 & 12 & 27&16\\
511$^2$ &16 & 9  & 12 & 12 & 36 
  &16\\
\hline
\end{tabular}
\end{tabular}
\caption{Number of iterations required by TGM and MGM - unilevel and two-level cases
(refer to (\ref{eq:def-rho}) for the definition of the constant $\rho$).} \label{table:1d-2d-multigrid}
\end{table}
\newline
Before dealing with other type of boundary conditions, we want to
give a  comparison of the performances of the proposed method with
respect to those achieved by considering, for instance, the
conjugate gradient (CG) method. Table \ref{table:cg-1D-2D}
reports, for increasing dimension, the Euclidean matrix condition
number $k_2(A_n+D_n)$, together with the number of iterations
required by the CG. As well known in the case $d0$, the CG method
requires all the $N(n)$ steps in order to reach the convergence.
Moreover, the non-structured part in the cases $d1$, $d2$, $d3$
increases the minimum eigenvalue of the resulting matrix so that
the whole condition number becomes moderate. As a consequence the
standard CG method is also effective. Notice that this good
behavior is no longer observed in the case $d4$, while our MGM
technique is still optimal. The same trend is observed in the
two-level setting.
%
\begin{table}[t!]
\footnotesize
\centering
{$B_n=\mathrm{tridiag}_n[-1, 2,-1]+$Diagonal}\\ 
\begin{tabular}[c]{|c| cc| cc| cc| cc| cc|}
\hline \multicolumn{1}{|c|}{$N(n)$} & \multicolumn{2}{|c|}{d0}
                          & \multicolumn{2}{|c|}{d1}
                          & \multicolumn{2}{|c|}{d2}
                          & \multicolumn{2}{|c|}{d3}
                          & \multicolumn{2}{|c|}{d4}\\
  & $k_2$ & nit & $k_2$ & nit & $k_2$ & nit & $k_2$ & nit & $k_2$ & nit\\
\hline
    31   & 4.14e+2 &  31  &  5.62e+0 & 18   &  7.98e+0 &22  &  8.16e+0 & 22  &  2.03e+1 & 27   \\
    63   & 1.65e+3 &  63  &  5.67e+0 & 18   &  8.02e+0 &21  &  8.19e+0 & 22  &  3.30e+1 & 34   \\
   127   & 6.63e+3 &  127 &  5.68e+0 & 18   &  8.02e+0 &21  &  8.19e+0 & 21  &  5.31e+1 & 43   \\
   255   & 2.65e+4 &  255 &  5.69e+0 & 17   &  8.02e+0 &21  &  8.20e+0 & 21  &  8.51e+1 & 54   \\
   511   & 1.06e+5 &  511 &  5.69e+0 & 17   &  8.05e+0 &21  &  8.20e+0 & 21  &  1.35e+2 & 66   \\
\hline
\end{tabular} \ \\ \ \\ \ \\
%
%
\centering {$B_n=\mathrm{tridiag}_{n^{(1)}}[-1, 2,
-1]\otimes I_{n^{(2)}}+I_{n^{(1)}} \otimes
\mathrm{tridiag}_{n^{(2)}}[-1, 2, -1] +$Diagonal}\\ 
\begin{tabular}[c]{|c| cc| cc| cc| cc| cc|}
\hline \multicolumn{1}{|c|}{$N(n)$} & \multicolumn{2}{|c|}{d0}
                          & \multicolumn{2}{|c|}{d1}
                          & \multicolumn{2}{|c|}{d2}
                          & \multicolumn{2}{|c|}{d3}
                          & \multicolumn{2}{|c|}{d4}\\
  & $k_2$ & nit & $k_2$ & nit & $k_2$ & nit & $k_2$ & nit & $k_2$ & nit\\
\hline
    31$^2$  & 4.14e+2 & 82   & 5.62e+0 & 18 & 7.98e+0 & 22 & 8.16e+0 & 22 & 3.83e+1 & 46 \\
    63$^2$  & 1.65e+3 & 163  & 5.67e+0 & 17 & 8.02e+0 & 22 & 8.19e+0 & 22 & 6.29e+1 & 57 \\
    127$^2$ & 6.63e+3 & 319  & 5.68e+0 & 16 & 8.02e+0 & 21 & 8.19e+0 & 21 & 1.01e+2 & 71 \\
    255$^2$ & 2.65e+4 & 623  & 5.69e+0 & 16 & 8.02e+0 & 21 & 8.20e+0 & 21 & 1.60e+2 & 83 \\
    511$^2$ & 1.06e+5 & 1215 & 5.69e+0 & 15 & 8.05e+0 & 21 & 8.20e+0 & 21 & 2.54e+2 & 99 \\
\hline
\end{tabular}
\caption{Euclidean condition number $k_2(A_n+D_n)$ and number of
iterations required by CG  - unilevel and two-level cases.}
\label{table:cg-1D-2D}
\end{table}
\subsubsection{Randomly generated indefinite corrections}
As a further interesting case, we want to test our proposal in the
case of randomly generated matrix corrections. More specifically,
we consider diagonal, symmetric tridiagonal, symmetric
pentadiagonal matrix corrections with random entries uniformly
distributed on the unit interval (cases $d5$, $d7$, and $d9$,
respectively) or normally distributed with mean zero and standard
deviation one (cases $d6$, $d8$, and $d10$, respectively). Notice
that in such a way we are also considering indefinite corrections.
Thus, in order to obtain a positive definite matrix $B_n$ and in
order to satisfy the crucial relation $A_n\le \vartheta B_n$ for
some positive $\vartheta$ independent of $n$, the arising random
matrices corrections are suitable scaled by $1/(\gamma n^2)$ in
the unilevel setting, with $\gamma$ being the number of non-zero
diagonals, and by $1/(\gamma (n^{(1)})^2)$ in the two-level
setting (assuming $n^{(1)}=n^{(2)}$).\newline Table
\ref{table:rand-randn-1D-2D} reports the Euclidean condition
number and the mean of the number of iterations required by the
MGM in the unilevel and two-level setting by considering, for each
case, ten examples of random matrix corrections.\newline All these
results confirm the effectiveness of our proposal. Though the
Euclidean condition numbers are fully comparable with those of the
$d0$ case, the number of required iterations does not worsen.
Conversely, the CG method requires for instance in the $d5$ case
$N(n)$ iterations in the unilevel setting, and $83, 163, 318, 621,
1212$ in the two-level one.\newline It is worth stressing that the
pentadiagonal corrections are reduced at the first projection to
tridiagonal matrices. More in general, bigger patterns are reduced
after few steps to a fixed pattern driven by the projector pattern
(see \cite{FS1,AD,ADS}).
\begin{table}[t!]
\footnotesize
\centering {$B_n=\mathrm{tridiag}_n[-1, 2, -1]+$random
correction}\\
\begin{tabular}[c]{|c| cc| cc| cc| cc| cc| cc|}
\hline \multicolumn{1}{|c|}{$N(n)$} & \multicolumn{2}{|c|}{d5}
                          & \multicolumn{2}{|c|}{d6}
                          & \multicolumn{2}{|c|}{d7}
                          & \multicolumn{2}{|c|}{d8}
                          & \multicolumn{2}{|c|}{d9}
                          & \multicolumn{2}{|c|}{d10}
                                                   \\
  & $k_2$ & \!\!\!\!\!\!nit\!\!\!\!\!\! & $k_2$ & \!\!\!\!\!\!nit\!\!\!\!\!\! & $k_2$ & \!\!\!\!\!\!nit\!\!\!\!\!\! & $k_2$
  & \!\!\!\!\!\!nit\!\!\!\!\!\!  & $k_2$ & \!\!\!\!\!\!nit\!\!\!\!\!\!  & $k_2$ & \!\!\!\!\!\!nit\!\!\!\!\!\!  \\
\hline
    31  &3.89e+2 & 3& 4.09e+2 & 3.5 & 3.83e+2 & 3 & 4.16e+2 & 3 & 3.84e+2 & 3 & 4.11e+2 & 3  \\
    63  &1.57e+3 & 7& 1.64e+3 & 7   & 1.53e+3 & 7 & 1.66e+3 & 7 & 1.51e+3 & 7 & 1.64e+3 & 7  \\
    12  &6.29e+3 & 8& 6.65e+3 & 8   & 6.18e+3 & 8 & 6.54e+3 & 8 & 6.11e+3 & 8 & 6.75e+3 & 8  \\
    255 &2.52e+4 & 8& 2.65e+4 & 8   & 2.46e+4 & 8 & 2.63e+4 & 8 & 2.44e+4 & 8 & 2.65e+4 & 8\\
    511 &1.01e+5 & 8& 1.05e+5 & 8   & 9.86e+4 & 8 & 1.06e+5 & 8 & 9.77e+4 & 8 & 1.05e+5 & 8 \\
\hline
\end{tabular} \ \\ \ \\ \ \\
%
%
\centering {$B_n=\mathrm{tridiag}_{n^{(1)}}[-1, 2,
-1]\otimes I_{n^{(2)}}+I_{n^{(1)}} \otimes
\mathrm{tridiag}_{n^{(2)}}[-1, 2, -1] +$random correction}\\
\begin{tabular}[c]{|c| cc| cc| cc| cc| cc| cc| }
\hline \multicolumn{1}{|c|}{$N(n)$} & \multicolumn{2}{|c|}{d5}
                          & \multicolumn{2}{|c|}{d6}
                          & \multicolumn{2}{|c|}{d7}
                          & \multicolumn{2}{|c|}{d8}
                          & \multicolumn{2}{|c|}{d9}
                          & \multicolumn{2}{|c|}{d10}
                                                  \\
  & $k_2$ & \!\!\!\!\!\!nit\!\!\!\!\!\! & $k_2$ & \!\!\!\!\!\!nit\!\!\!\!\!\! & $k_2$ & \!\!\!\!\!\!nit\!\!\!\!\!\! & $k_2$
  & \!\!\!\!\!\!nit\!\!\!\!\!\!  & $k_2$ & \!\!\!\!\!\!nit\!\!\!\!\!\!  & $k_2$ & \!\!\!\!\!\!nit\!\!\!\!\!\!  \\
\hline
    31$^2$  & 4.02e+2 & 16 & 4.15e+2 & 16 & 3.98e+2 & 16 & 4.15e+2 & 16 & 3.96e+2 & 16 & 4.14e+2 & 16   \\
    63$^2$  & 1.61e+3 & 16 & 1.66e+3 & 16 & 1.59e+3 & 16 & 1.66e+3 & 16 & 1.59e+3 & 16 & 1.65e+3 & 16   \\
    127$^2$ & 6.47e+3 & 16 & 6.64e+3 & 16 & 6.39e+3 & 16 & 6.63e+3 & 16 & 6.36e+3 & 16 & 6.63e+3 & 16  \\
    255$^2$ & 2.58e+4 & 16 & 2.65e+4 & 16 & 2.55e+4 & 16 & 2.65e+4 & 16 & 2.54e+4 & 16 & 2.65e+4 & 16   \\
    511$^2$ & 1.03e+5 & 16 & 1.06e+5 & 16 & 1.02e+5 & 16 & 1.06e+5 & 16 & 1.01e+5 & 16 & 1.06e+5 & 16  \\
\hline
\end{tabular}
\caption{Euclidean condition number $k_2(A_n+D_n)$ and mean number
of iterations required by MGM - unilevel and two-level cases.}
\label{table:rand-randn-1D-2D}
\end{table}
%
\subsubsection{Periodic and Reflective boundary conditions}
In the case of periodic boundary conditions the obtained matrix
sequence is the unilevel circulant matrix sequence $\{S_n(f)\}_n$ generated
by the function $f(t)=2-2\cos(t)$, $t\in (0,2\pi]$. Following
\cite{mcirco}, we consider the operator $T_0^1\in
\mathbb{R}^{n_0\times n_1}$, $n_0=2n_1$, such that
\[
  (T_0^1)_{i,j}=\left\{ \begin{array}{ll}
            1 & \ \ \mbox{\textrm for }\  i=2j-1, \ \ j=1,\ldots,n_1,  \\
            0 & \ \ \mbox{\textrm otherwise},
          \end{array} \right.
\]
and we define a projector $(p_0^1)^H$, $p_0^1\in
\mathbb{R}^{n_0\times n_1}$, as
$p_0^1= P_0 T_0^1$,  $P_0=\ {S}_0(p)$, $p(t)=2+2\cos(t)$.
It must be outlined that in the $d0$ case the arising matrices are
singular, so that we consider the classical Strang correction
\cite{T-LAA-1995}
\[
\tilde S_{n_0}(f)= S_{n_0}(f)+ f\left(\frac{2\pi}{N({n_0})}\right)
\frac{ee^t}{N({n_0})},
\]
where $e$ is the vector of all ones.
By using tensor arguments, our results plainly
extend to the two-level setting
and the results in the top part of Table \ref{table:2dmultigrid-bis} confirm the optimality of the
proposed TGM and its extension to MGM (the case $d4$ requires to set
$\rho=1$).
\newline
When dealing with reflective boundary conditions, the obtained
matrix sequence is the unilevel  DCT III matrix sequence ${C_n(f)}_n$
generated by the function $f(t)=2-2\cos(t)$, $t\in (0,2\pi]$.
Following \cite{mcoseni}, we consider the operator $T_0^1\in
\mathbb{R}^{n_0\times n_1}$, $n_0=2n_1$, such that
\[
  (T_0^1)_{i,j}=\left\{ \begin{array}{ll}
            1 & \ \ \mbox{\textrm for }\  i\in \{2j-1,2j\}, \ \ j=1,\ldots,n_1,  \\
            0 & \ \ \mbox{\textrm otherwise},
          \end{array} \right.
\]
and we define a projector $(p_0^1)^H$, $p_0^1\in
\mathbb{R}^{n_0\times n_1}$, as
$p_0^1= P_0 T_0^1$, $P_0=C_0(p)$, $p(t)=2+2\cos(t)$.
Again, the two-level setting is treated by using tensor arguments.
The results in bottom part of Table \ref{table:2dmultigrid-bis} confirm the optimality of the proposed
TGM and its extension to MGM. It is worth stressing that in the
$d0$ case we are considering the matrix
\[
\tilde C_{n_0}(f)=C_{n_0}(f)+f\left(\frac{\pi}{N({n_0})}\right)
\frac{ee^t}{N({n_0})}.
\]
The  numerical tests related to the unilevel Laplacian with
periodic or reflective boundary conditions have the same flavor.
\begin{table}[t!]
\footnotesize
\centering
\begin{tabular}[c]{ccc} 
\multicolumn{3}{c}{{$B_n=$ two-level circulant $S_n(f)+$
Diagonal, $f(t_1,t_2)=4-2\cos(t_1)-2\cos(t_2)$}}\\  
\begin{tabular}[c]{|c|c|c|c|c|c| }\hline
\multicolumn{6}{|c|}{TGM} \\ \hline
\multicolumn{1}{|c|}{$N(n)$} & \multicolumn{1}{|c|}{d0}
                          & \multicolumn{1}{|c|}{d1}
                          & \multicolumn{1}{|c|}{d2}
                          & \multicolumn{1}{|c|}{d3}
                          & \multicolumn{1}{|c|}{d4}
\\
\multicolumn{1}{|c|}{} & \multicolumn{1}{|c|}{}
                          & \multicolumn{1}{|c|}{}
                          & \multicolumn{1}{|c|}{}
                          & \multicolumn{1}{|c|}{}
                          & \multicolumn{1}{|c|}{}
                          \\
\hline
32$^2$   & 15 & 8 & 11 & 11 & 14 \\
64$^2$   & 15 & 7 & 11 & 11 & 15 \\
128$^2$  & 15 & 7 & 11 & 11 & 15 \\
256$^2$  & 15 & 7 & 11 & 11 & 15 \\
512$^2$  & 15 & 7 & 11 & 11 & 15 \\
\hline
\end{tabular}
& \ &
\begin{tabular}[c]{|c|c|c|c|c|c|c|}\hline
\multicolumn{7}{|c|}{MGM} \\ \hline
\multicolumn{1}{|c|}{$N(n)$} & \multicolumn{1}{|c|}{d0}
                          & \multicolumn{1}{|c|}{d1}
                          & \multicolumn{1}{|c|}{d2}
                          & \multicolumn{1}{|c|}{d3}
                          & \multicolumn{1}{|c|}{d4} &
\\
\multicolumn{1}{|c|}{} & \multicolumn{1}{|c|}{}
                          & \multicolumn{1}{|c|}{}
                          & \multicolumn{1}{|c|}{}
                          & \multicolumn{1}{|c|}{}
                          & \multicolumn{1}{|c|}{$\rho\!=\!0$} &  \multicolumn{1}{|c|}{$\rho\!=\!1$}
                          \\
\hline
32$^2$  & 15  & 8 & 11 & 11 & 14&14\\
64$^2$  & 15  & 7 & 11 & 11 & 15&15\\
128$^2$ & 15  & 7 & 11 & 11 & 16&14\\
256$^2$ & 15  & 7 & 11 & 11 & 24&14\\
512$^2$ & 15  & 7 & 11 & 11 & 34&14\\
\hline
\end{tabular}
\end{tabular} \ \\ \ \\
\begin{tabular}[c]{ccc} 
\multicolumn{3}{c}{{$B_n=$ two-level DCT III $C_n(f)+$
Diagonal, $f(t_1,t_2)=4-2\cos(t_1)-2\cos(t_2)$}}\\ 
\begin{tabular}[c]{|c|c|c|c|c|c|}\hline
\multicolumn{6}{|c|}{TGM}\\ \hline
\multicolumn{1}{|c|}{$N(n)$} & \multicolumn{1}{|c|}{d0}
                          & \multicolumn{1}{|c|}{d1}
                          & \multicolumn{1}{|c|}{d2}
                          & \multicolumn{1}{|c|}{d3}
                          & \multicolumn{1}{|c|}{d4}
\\
\multicolumn{1}{|c|}{} & \multicolumn{1}{|c|}{}
                          & \multicolumn{1}{|c|}{}
                          & \multicolumn{1}{|c|}{}
                          & \multicolumn{1}{|c|}{}
                          & \multicolumn{1}{|c|}{}
                          \\
\hline
32$^2$   & 16 & 6 & 10 & 10 & 12 \\
64$^2$   & 16 & 6 & 10 & 10 & 11 \\
128$^2$  & 16 & 5 & 10 & 10 & 11 \\
256$^2$  & 16 & 5 & 9  & 9  & 11 \\
512$^2$  & 16 & 5 & 9  & 9  & 11 \\
\hline
\end{tabular}
& \ &
\begin{tabular}[c]{|c|c|c|c|c|c| c|c|c|c|c|c|c|}\hline
\multicolumn{7}{|c|}{MGM} \\ \hline
\multicolumn{1}{|c|}{$N(n)$} & \multicolumn{1}{|c|}{d0}
                          & \multicolumn{1}{|c|}{d1}
                          & \multicolumn{1}{|c|}{d2}
                          & \multicolumn{1}{|c|}{d3}
                          & \multicolumn{2}{|c|}{d4}
\\
\multicolumn{1}{|c|}{} & \multicolumn{1}{|c|}{}
                          & \multicolumn{1}{|c|}{}
                          & \multicolumn{1}{|c|}{}
                          & \multicolumn{1}{|c|}{}
                          & \multicolumn{1}{|c|}{$\rho\!=\!0$} &  \multicolumn{1}{|c|}{$\rho\!=\!1$}
                          \\
\hline
32$^2$  & 16  & 6 & 10 & 10 & 12 &12\\
64$^2$  & 16  & 6 & 10 & 10 & 11 &11\\
128$^2$ & 16  & 5 & 10 & 10 & 11 &10\\
256$^2$ & 16  & 5 & 9  & 9  & 17 &9\\
512$^2$ & 16  & 5 & 9  & 9  & 27 &9\\
\hline
\end{tabular}
\end{tabular} 
\caption{Number of iterations required by TGM and MGM - two-level
cases (refer to (\ref{eq:def-rho}) for the definition of the
constant $\rho$).} \label{table:2dmultigrid-bis}
\end{table}
%
\subsection{Other examples}
In this section we give numerical evidences of the optimality of
TGM and MGM results in a more general setting.
\subsubsection{Higher order $\tau$ discretizations plus diagonal systems}
We consider $\tau$ matrix sequences arising from the
discretization of higher order differential problems with proper
homogeneous boundary conditions on $\partial\Omega$:
\begin{equation} \label{eq:modkb2D}
(-1)^q \sum_{i=1}^d {\partial^{2q}\over \partial
x_i^{2q}}u(x)+\mu(x)u(x)=h(x) \ \mbox{on}\ \Omega=(0,1)^d,
\end{equation}
i.e., $\{B_n=A_n+D_n\}_n$, where $A_n=\tau_n(f)$ with
$f(t)=\sum_{i=1}^d (2-2\cos(t_i))^q$. More specifically, in the
unilevel case we define $p(t)=[2+2\cos(t)]^w$ where $w$ is chosen
according to conditions in \cite{FS1,CCS,Sun}:
in order to have a MGM optimality we must take $w$ at least equal to
$1$ if $q=1$ and at least equal to $2$ if $q=2,3$.
\newline With respect to the two-level
problem, we consider again the most trivial extension (and less
expensive from a computational point of view) of the unilevel
projector to the two-level setting, given by  $P_n=\tau_n(p)$ with
$p(t_1,t_2)=[(2+2\cos(t_1))(2+2\cos(t_2))]^w$, $w=1,2,3$. \newline
Clearly, the
lower is the value of $w$, the greater will be the advantage from a
computational viewpoint. Indeed, Table \ref{table:multigrid:tau-ho}
confirms the need of these constraints with respect to the case
$q=2$, this being the only $d0$ case where we observe a growth in
the number of iterations with respect to $N(n)$. Nevertheless, it
should be noticed that in the same case the contribution of the
non-structured part improves the numerical behavior since the
minimal eigenvalue is increased.
\newline The remaining results in Table \ref{table:multigrid:tau-ho}
confirm the optimality of the corresponding MGM (the $d4$ case
requires to set $\rho$ in a proper way as just observed in the
Laplacian case). Notice that the bandwidth of the
non-structured diagonal correction is increased by subsequent
projections until a maximal value corresponding to $4w-1$ is
reached (for a discussion on the evolution of the bandwidth when a
generic (multilevel) band system is encountered see
\cite{FS1,AD,ADS}).
\begin{table}[t!]
\footnotesize
\centering
{\centering {$B_n=$two-level $\tau+$ Diagonal,
$f(t_1,t_2)=(2-2\cos(t_1))^q+(2-2\cos(t_2))^q$}\\ }
\begin{tabular}[c]{|c|c|c|c|c|c|c| c|c|c|c|c|c|c|}\hline
\multicolumn{14}{|c|}{$q=2$}  \\
\hline \multicolumn{7}{|c|}{$w=1$} & \multicolumn{7}{|c|}{$w=2$}
\\ \hline
\multicolumn{1}{|c|}{$N(n)$} & \multicolumn{1}{|c|}{d0}
                          & \multicolumn{1}{|c|}{d1}
                          & \multicolumn{1}{|c|}{d2}
                          & \multicolumn{1}{|c|}{d3}
                          & \multicolumn{2}{|c|}{d4} &
 \multicolumn{1}{|c|}{$N(n)$} &
\multicolumn{1}{|c|}{d0}
                          & \multicolumn{1}{|c|}{d1}
                          & \multicolumn{1}{|c|}{d2}
                          & \multicolumn{1}{|c|}{d3}
                          & \multicolumn{2}{|c|}{d4}
\\
\multicolumn{1}{|c|}{} & \multicolumn{1}{|c|}{}
                          & \multicolumn{1}{|c|}{}
                          & \multicolumn{1}{|c|}{}
                          & \multicolumn{1}{|c|}{}
                          & \multicolumn{1}{|c|}{$\rho\!=\!0$} &  \multicolumn{1}{|c|}{$\rho\!=\!2$}
                          &
\multicolumn{1}{|c|}{} & \multicolumn{1}{|c|}{}
                          & \multicolumn{1}{|c|}{}
                          & \multicolumn{1}{|c|}{}
                          & \multicolumn{1}{|c|}{}
                          & \multicolumn{1}{|c|}{$\rho\!=\!0$} &  \multicolumn{1}{|c|}{$\rho\!=\!2$}
                          \\
\hline
31$^2$   & 37 & 29 & 32 & 32 & 36 &36
& 31$^2$ & 35 & 28 & 31 & 31 & 34 &34
\\
63$^2$   & 44 & 28 & 32 & 32 & 43 &35
& 63$^2$ & 36 & 28 & 31 & 31 & 34 &34
\\
127$^2$  & 80 & 28 & 31 & 31 & 73 &35
& 127$^2$& 36 & 27 & 30 & 30 & 56 &34
\\
255$^2$  & 140& 27 & 30 & 30 & 109&35
& 255$^2$ & 36 & 27 & 30 & 30 & 89 &33
\\
511$^2$  & 235 & 27 & 30 & 30 &  151  &35
& 511$^2$ & 36 & 27 & 30 & 30 & 129   &33\\
\hline
\multicolumn{14}{|c|}{$q=3$}  \\
\hline \multicolumn{7}{|c|}{$w=2$} & \multicolumn{7}{|c|}{$w=3$} \\
\hline \multicolumn{1}{|c|}{$N(n)$} & \multicolumn{1}{|c|}{d0}
                          & \multicolumn{1}{|c|}{d1}
                          & \multicolumn{1}{|c|}{d2}
                          & \multicolumn{1}{|c|}{d3}
                          & \multicolumn{2}{|c|}{d4} &
 \multicolumn{1}{|c|}{$N(n)$} & \multicolumn{1}{|c|}{d0}
                          & \multicolumn{1}{|c|}{d1}
                          & \multicolumn{1}{|c|}{d2}
                          & \multicolumn{1}{|c|}{d3}
                          & \multicolumn{2}{|c|}{d4}
\\
\multicolumn{1}{|c|}{} & \multicolumn{1}{|c|}{}
                          & \multicolumn{1}{|c|}{}
                          & \multicolumn{1}{|c|}{}
                          & \multicolumn{1}{|c|}{}
                          & \multicolumn{1}{|c|}{$\rho\!=\!0$} &  \multicolumn{1}{|c|}{$\rho\!=\!2$}
                          &
\multicolumn{1}{|c|}{} & \multicolumn{1}{|c|}{}
                          & \multicolumn{1}{|c|}{}
                          & \multicolumn{1}{|c|}{}
                          & \multicolumn{1}{|c|}{}
                          & \multicolumn{1}{|c|}{$\rho\!=\!0$} &  \multicolumn{1}{|c|}{$\rho\!=\!2$}
                          \\
\hline
31$^2$   & 72 & 64 & 67 & 67 & 69  &69 
& 31$^2$ & 68 & 61 & 64 & 64 & 66 &66 
\\
63$^2$   & 73 & 65 & 68 & 68 & 71  &71
& 63$^2$ & 72 & 64 & 67 & 67 & 70  &70
\\
127$^2$  & 73 & 64 & 68 & 68 & 137 &71
& 127$^2$& 72 & 63 & 67 & 67 & 127 &70
\\
255$^2$  & 73 & 64 & 67 & 67 & 193 &70
& 255$^2$& 72 & 63 & 67 & 67 & 182 &70
\\
511$^2$  & 73 & 64 & 67 & 67 & 276 &70          & 511$^2$ & 72 & 63 & 67 & 67 & 260    &70 \\
\hline
\end{tabular}
\caption{Number of required MGM iterations - 
two-level cases (refer to (\ref{eq:def-rho}) for the definition of
the constant $\rho$).} \label{table:multigrid:tau-ho}
\end{table}
%
\subsubsection{Higher order circulant discretizations plus diagonal systems}
We consider circulant matrix sequences arising from the
approximation of higher order differential problems with proper
homogeneous periodic boundary conditions  on $\partial\Omega$ as in
(\ref{eq:modkb2D}), i.e., $\{B_n=A_n+D_n\}_n$, where $A_n=S_n(f)$
with $f(t)=\sum_{i=1}^d (2-2\cos(t_i))^q$. The choice of the
generating function for the projector is the same as in the previous
section (see \cite{mcirco}). Indeed, Table
\ref{table:multigrid:circulant-ho} shows the importance of these
constraints with respect to the case $d0$ with $q=2$. It is worth
mentioning that the optimality of the corresponding MGM is again
confirmed (for the case $d4$ the parameter $\rho$ has to be set in a
proper way).
%
\begin{table}[t!]
\footnotesize
\centering
{\centering {$B_n=$ two-level circulant $+$ Diagonal,
$f(t_1,t_2)=(2-2\cos(t_1))^q+(2-2\cos(t_2))^q$}\\ } 
\begin{tabular}[c]{|c|c|c|c|c|c|c| c|c|c|c|c|c|c|}\hline
\multicolumn{14}{|c|}{$q=2$}  \\
\hline \multicolumn{7}{|c|}{$w=1$} & \multicolumn{7}{|c|}{$w=2$} \\
\hline \multicolumn{1}{|c|}{$N(n)$} & \multicolumn{1}{|c|}{d0}
                          & \multicolumn{1}{|c|}{d1}
                          & \multicolumn{1}{|c|}{d2}
                          & \multicolumn{1}{|c|}{d3}
                          & \multicolumn{2}{|c|}{d4} &
 \multicolumn{1}{|c|}{$N(n)$} &
\multicolumn{1}{|c|}{d0}
                          & \multicolumn{1}{|c|}{d1}
                          & \multicolumn{1}{|c|}{d2}
                          & \multicolumn{1}{|c|}{d3}
                          & \multicolumn{2}{|c|}{d4}
\\
\multicolumn{1}{|c|}{} & \multicolumn{1}{|c|}{}
                          & \multicolumn{1}{|c|}{}
                          & \multicolumn{1}{|c|}{}
                          & \multicolumn{1}{|c|}{}
                          & \multicolumn{1}{|c|}{$\rho\!=\!0$} &  \multicolumn{1}{|c|}{$\rho\!=\!2$}
                          &
\multicolumn{1}{|c|}{} & \multicolumn{1}{|c|}{}
                          & \multicolumn{1}{|c|}{}
                          & \multicolumn{1}{|c|}{}
                          & \multicolumn{1}{|c|}{}
                          & \multicolumn{1}{|c|}{$\rho\!=\!0$} &  \multicolumn{1}{|c|}{$\rho\!=\!1$}
                          \\
\hline
32$^2$   & 34  & 23 & 26 & 27  & 32& 32
& 32$^2$   & 34 & 23 & 27 & 27 & 31 & 31\\
64$^2$   & 42  & 22 & 25 & 26  & 38& 33
& 64$^2$   & 33 & 22 & 25 & 25 & 31 & 31\\
128$^2$  & 76  & 22 & 25 & 25  & 56& 33
& 128$^2$  & 33 & 22 & 25 & 25 & 44 & 31\\
256$^2$  & 130 & 21 & 24 & 24  & 75& 34
& 256$^2$  & 33 & 21 & 24 & 24 & 62 & 31\\
512$^2$  & 213 & 21 & 25 & 225 & 96  & 34
& 512$^2$  & 33 & 20 & 24 & 24 & 85 & 31\\
\hline
\multicolumn{14}{|c|}{$q=3$}  \\
\hline \multicolumn{7}{|c|}{$w=2$} & \multicolumn{7}{|c|}{$w=3$} \\
 \hline
\multicolumn{1}{|c|}{$N(n)$} & \multicolumn{1}{|c|}{d0}
                          & \multicolumn{1}{|c|}{d1}
                          & \multicolumn{1}{|c|}{d2}
                          & \multicolumn{1}{|c|}{d3}
                          & \multicolumn{2}{|c|}{d4} &
 \multicolumn{1}{|c|}{$N(n)$} &
\multicolumn{1}{|c|}{d0}
                          & \multicolumn{1}{|c|}{d1}
                          & \multicolumn{1}{|c|}{d2}
                          & \multicolumn{1}{|c|}{d3}
                          & \multicolumn{2}{|c|}{d4}
\\
\multicolumn{1}{|c|}{} & \multicolumn{1}{|c|}{}
                          & \multicolumn{1}{|c|}{}
                          & \multicolumn{1}{|c|}{}
                          & \multicolumn{1}{|c|}{}
                          & \multicolumn{1}{|c|}{$\rho\!=\!0$} &  \multicolumn{1}{|c|}{$\rho\!=\!1$}
                          &
\multicolumn{1}{|c|}{} & \multicolumn{1}{|c|}{}
                          & \multicolumn{1}{|c|}{}
                          & \multicolumn{1}{|c|}{}
                          & \multicolumn{1}{|c|}{}
                          & \multicolumn{1}{|c|}{$\rho\!=\!0$} &  \multicolumn{1}{|c|}{$\rho\!=\!1$}
                          \\
\hline
32$^2$   & 68  & 56 & 59 & 60 & 65  & 65& 32$^2$   & 68 & 56 & 59 & 60 & 65  &65\\
64$^2$   & 67  & 55 & 58 & 58 & 64  & 64& 64$^2$   & 66 & 54 & 58 & 58 & 64  &63\\
128$^2$  & 67  & 54 & 58 & 58 & 88  & 64& 128$^2$  & 66 & 53 & 57 & 58 & 82  &63\\
256$^2$  & 67  & 53 & 57 & 57 & 118 & 65& 256$^2$  & 66 & 52 & 57 & 57 & 111 &63\\
512$^2$  & 67  & 52 & 57 & 57 & 158 & 65& 512$^2$  & 66 & 51 & 56 & 57 & 149 &63\\
\hline
\end{tabular}
\caption{Number of required MGM iterations - 
two-level cases (refer to (\ref{eq:def-rho}) for the definition of
the constant $\rho$).}  \label{table:multigrid:circulant-ho}
\end{table}
\subsubsection{Reflective BCs discretizations plus diagonal
systems}\label{subsez:zeropi-cosine}
We consider an example of DCT-III matrix sequences arising from
the discretization of integral problems with reflective boundary
conditions (see \cite{NCT}), i.e., $\{B_n=A_n+D_n\}_n$, where
$A_n=C_n(f)$ with $f$ having nonnegative Fourier coefficients as
it is required for the point spread function in the modelling of
image blurring, see \cite{BB}. A simple model is represented by
$f(t)=f_d(t):=\sum_{i=1}^d (2+2\cos(t_i))$ where, by the way, the
product $f_1(t_1)f_1(t_2)$ is encountered when treating
super-resolution or high resolution problems, see e.g.
\cite{super}. The choice of the generating function for the
projector is the same as in \cite{mcoseni}. \newline The results
in Table \ref{table:multigrid:cosined2zeropi} confirm again the
optimality of the corresponding MGM.
%
\begin{table}[t!]
\footnotesize \centering
{\centering {$B_n=$ two-level DCT III $C_n(f)+$
Diagonal,
$f(t_1,t_2)=(2+2\cos(t_1))+(2-2\cos(t_2))$.} \\ } 
%
\begin{tabular}[c]{ccc}
\begin{tabular}[c]{|c|c|c|c|c|c| c|c|c|c|c|c|}\hline
\multicolumn{6}{|c|}{TGM} \\
\hline
\multicolumn{1}{|c|}{$N(n)$} & \multicolumn{1}{|c|}{d0}
                          & \multicolumn{1}{|c|}{d1}
                          & \multicolumn{1}{|c|}{d2}
                          & \multicolumn{1}{|c|}{d3}
                          & \multicolumn{1}{|c|}{d4}
\\
\hline
32$^2$   & 5 & 6 & 7 & 8 & 6 \\
64$^2$   & 4 & 5 & 7 & 8 & 6  \\
128$^2$  & 4 & 5 & 7 & 8 & 6  \\
256$^2$  & 4 & 5 & 7 & 8 & 6  \\
512$^2$  & 4 & 5 & 7 & 8 & 6  \\
\hline
\end{tabular}
& \ &
%
\begin{tabular}[c]{|c|c|c|c|c|c| c|c|c|c|c|c|}\hline
\multicolumn{6}{|c|}{MGM} \\
\hline
\multicolumn{1}{|c|}{$N(n)$} & \multicolumn{1}{|c|}{d0}
                          & \multicolumn{1}{|c|}{d1}
                          & \multicolumn{1}{|c|}{d2}
                          & \multicolumn{1}{|c|}{d3}
                          & \multicolumn{1}{|c|}{d4}
\\
\hline
32$^2$  & 5 & 6 & 7 & 8 & 6 \\
64$^2$  & 5 & 5 & 9 & 9 & 6 \\
128$^2$ & 4 & 5 & 9 & 9 & 6 \\
256$^2$ & 4 & 5 & 9 & 9 & 6 \\
512$^2$ & 4 & 5 & 9 & 9 & 6 \\
\hline
\end{tabular}
\end{tabular}
\caption{Number of iterations required by TGM and MGM - unilevel
and two level cases (refer to (\ref{eq:def-rho}) for the
definition of the constant $\rho$).}
\label{table:multigrid:cosined2zeropi}
\end{table}
\section{Concluding Remarks and Perspectives}\label{sec:finrem}
The algebraic tools given in Section \ref{sec:extension} and
Section \ref{sec:laplacian} revealed that, if a suitable TGM for a
Hermitian positive definite matrix sequence $\{A_n\}_n$ is
available and another Hermitian positive definite uniformly
bounded sequence $\{B_n\}_n$ is given such that $A_n \le \vartheta
B_n$, for $n$ large enough, then the same strategy works almost
unchanged for $\{B_n\}_n$ too. As an example, this means that if
the method is optimal for the first sequence then it is optimal
for the second as well. The same results should hold for the MGM
procedures, but here only a wide set of numerical evidences has
been provided for supporting this claim: the related theory will be
a subject of future investigations taking into account the final
remarks in Section \ref{tgm-opt} and the discussion in Section
\ref{mgm-opt}.
\newline We point out that the latter goal is quite important.
Indeed, it is not difficult to prove relations of the form
$\vartheta_1 A_n\le B_n\le  \vartheta_2 A_n$ with $B_n$ being
discretization of an elliptic variable coefficient problem, $A_n$
being the same discretization in the constant coefficient case, and
where $\vartheta_1,\vartheta_2$ are positive constants independent
of $n$ and mainly depending on the ellipticity parameters of the
problem. Therefore, the above mentioned results would represent a
link for inferring MGM optimality on a general (possibly high order)
variable coefficient elliptic problem, starting from the MGM
optimality for the structured part, i.e., the one related to the
constant coefficient discretization.
\newline Finally, we point out that the latter idea has been used
essentially for the structured plus diagonal systems coming from
approximated elliptic partial differential equations with different
boundary conditions. However, the same approach is applicable to a
wide variety of cases, as sketched for instance in Section
\ref{subsez:zeropi-cosine}.
\appendix
\section{Appendix}\label{sez:appendix}
For reader convenience, we report the essential steps of the proof
of Theorems  \ref{teo:TGMconv} and
\ref{teo:TGMconv-pre-post}. 
Let us start by proving Theorem \ref{teo:TGMconv}. As
demonstrated in Theorem 5.2 in \cite{RStub}, the existence of
$\beta>0$ such that $
\min_{y\in \mathbb{C}^{N(n_1)} } \| x -p_0^1 y \|_{D_0}^2 \le
\beta \ \| x \|_{A_0}^2$  for all $x\in \mathbb{C}^{N(n_0)}$
implies the validity of the so called \emph{approximation property} only in the range of $CGC_0$,
i.e.,  the existence of $\beta>0$ such that
\begin{equation}\label{eq:approximation-property}
\|CGC_0 x\|_{A_0}^2 \le \beta \ \|CGC_0 x\|_{A_0D_0^{-1}A_0}^2
\quad \forall x\in \mathbb{C}^{N(n_0)}
\end{equation}
where $CGC_0= I_0-p_0^1 A_1^{-1}(p_0^1)^H A_0$.
Thus, by virtue of the \emph{post-smoothing property}
(\ref{hp:1post-0}) and of (\ref{eq:approximation-property}), for all $x\in \mathbb{C}^{N(n_0)}$ we find
\begin{eqnarray}
\|V_{0,\mathrm{post}}CGC_0 x\|_{A_0}^2 &\le& \|CGC_0
x\|_{A_0}^2-\alpha_\mathrm{post}\
\|CGC_0  x\|_{A_0D_0^{-1}A_0}^2 \nonumber\\
&\le& \|CGC_0  x\|_{A_0}^2-\frac{\alpha_\mathrm{post}}{\beta}\
\|CGC_0
 x\|_{A_0}^2 \nonumber \\
&=& \left (1-\frac{\alpha_\mathrm{post}}{\beta} \right)\|CGC_0
x\|_{A_0}^2 \nonumber \\
&\le & \left (1-\frac{\alpha_\mathrm{post}}{\beta}
\right)\|x\|_{A_0}^2\label{eq:primaparte-0},
\end{eqnarray}
being $\|CGC_0 \|_{A_0}=1$. 
Since $TGM_0=V_{0,\mathrm{post}}CGC_0$ in the case where no
pre-smoothing is considered, the latter is the same as $\|TGM_0\|_{A_0}
\le \sqrt{1-\frac{\alpha_\mathrm{post}}{\beta}}$ and hence Theorem
\ref{teo:TGMconv} is proved.
\par
Now let us prove Theorem \ref{teo:TGMconv-pre-post}. Since the \emph{approximation property}
(\ref{hp:2-1}) implies clearly (\ref{eq:approximation-property}),
by repeating the very same steps as before and exploiting the
\emph{post-smoothing property} (\ref{hp:1post-1}), for all $x\in
\mathbb{C}^{N(n_0)}$ we find
\begin{eqnarray}
\|V_{0,\mathrm{post}}CGC_0 V_{0,\mathrm{pre}} x\|_{A_0}^2 &\le&
\|CGC_0 V_{0,\mathrm{pre}} x\|_{A_0}^2-\alpha_\mathrm{post}\
\|CGC_0 V_{0,\mathrm{pre}} x\|_{A_0D_0^{-1}A_0}^2 \nonumber\\
&\le& \|CGC_0 V_{0,\mathrm{pre}}
x\|_{A_0}^2-\frac{\alpha_\mathrm{post}}{\beta}\ \|CGC_0
V_{0,\mathrm{pre}} x\|_{A_0}^2 \nonumber \\
&=& \left (1-\frac{\alpha_\mathrm{post}}{\beta} \right)\|CGC_0
V_{0,\mathrm{pre}} x\|_{A_0}^2. \label{eq:primaparte}
\end{eqnarray}
In addition, by using (\ref{hp:2-1}) and the
\emph{pre-smoothing property} (\ref{hp:1pre-1}), respectively, for all $x\in
\mathbb{C}^{N(n_0)}$ we obtain
\begin{eqnarray*}
\|CGC_0 V_{0,\mathrm{pre}} x\|_{A_0}^2 &\le& \beta  \|
V_{0,\mathrm{pre}} x\|_{A_0D_0^{-1}A_0}^2\\
\|V_{0,\mathrm{pre}} x\|_{A_0D_0^{-1}A_0}^2  &\le&
\alpha_{\mathrm{pre}}^{-1} ( \|x\|_{A_0}^2 -\| V_{0,\mathrm{pre}}
x\|_{A_0}^2).
\end{eqnarray*}
Hence
\begin{eqnarray*}
\frac{\alpha_{\mathrm{pre}}}{\beta}\|CGC_0 V_{0,\mathrm{pre}}
x\|_{A_0}^2 &\le& \|x\|_{A_0}^2 -\| V_{0,\mathrm{pre}}
x\|_{A_0}^2 
\le \|x\|_{A_0}^2 -\| CGC_0 V_{0,\mathrm{pre}} x\|_{A_0}^2,
\end{eqnarray*}
since 
$
\| CGC_0 V_{0,\mathrm{pre}} x\|_{A_0}^2 \le \| CGC_0\|_{A_0}^2\|
V_{0,\mathrm{pre}} x\|_{A_0}^2 =  \| V_{0,\mathrm{pre}} x\|_{A_0}^2,
$
being $\| CGC_0\|_{A_0}=1$.  Therefore, for all $x\in
\mathbb{C}^{N(n_0)}$, it holds
\begin{equation} \label{eq:secondaparte}
\|CGC_0 V_{0,\mathrm{pre}} x\|_{A_0}^2 \le
\left(1+\frac{\alpha_{\mathrm{pre}}}{\beta} \right)^{-1}
\|x\|_{A_0}^2.
\end{equation}
By using inequality (\ref{eq:secondaparte}) in (\ref{eq:primaparte}), we
have
\begin{equation*}
\|V_{0,\mathrm{post}}CGC_0 V_{0,\mathrm{pre}} x\|_{A_0}^2 \le
\frac{1-{\alpha_\mathrm{post}}/{\beta}}
{1+{\alpha_\mathrm{pre}}/{\beta}}
 \|x\|_{A_0}^2,
\end{equation*}
and the proof of Theorem \ref{teo:TGMconv-pre-post} is
concluded.
%
%


%

\begin{thebibliography}{99}
\bibitem{AD}
{A. Aric\`o and M. Donatelli}, \textit{A V-cycle multigrid for
multilevel matrix algebras: proof of optimality}, Numer. Math.,
{105} (2007),  no. 4, pp.~511--547.
%
\bibitem{ADS}
{A. Aric\`o, M. Donatelli, and S. Serra-Capizzano},
\textit{V-cycle optimal convergence for certain (multilevel)
structured linear system}, SIAM J. Matrix Anal. Appl., {26-1}
(2004), pp.~186--214.
%
\bibitem{Axelsson}
{O. Axelsson}, \textit{Iterative Solution Methods}, Cambridge
University Press, Cambridge, 1994.
%
\bibitem{AxN}
{O. Axelsson and M. Neytcheva}, \textit{The algebraic multilevel
iteration methods -- theory and applications}, {Proc. 2nd Int.
Coll. on Numerical Analysis}, D. Bainov Ed., Plovdiv (Bulgaria),
August 1993, pp.~13--23.
%
\bibitem{BGST} D. Bertaccini, G. H. Golub, S. Serra-Capizzano, and C. Tablino-Possio,
\textit{Preconditioned HSS methods for the solution of non-Hermitian
positive definite linear systems and applications to the discrete
convection-diffusion equation}, Numer. Math., 99 (2005), pp.
441--484.
%
\bibitem{BB}
{M. Bertero and P. Boccacci}, \textit{Introduction to Inverse
Problems in Imaging}, Inst. of Physics Publ. Bristol and
Philadelphia, London, 1998.
%
\bibitem{bhatia}
 R. Bhatia, \textit{Matrix Analysis},
 Springer Verlag, New York, 1997.
%
\bibitem{BC}
{D. Bini and M. Capovani}, \textit{Spectral and computational
properties of band symmetric Toeplitz matrices}, {Linear Algebra
Appl.}, {52/53} (1983), pp.~99--125.
%
\bibitem{CCS}
{R.~H. Chan, Q. Chang, and H. Sun}, \textit{Multigrid method for
ill-conditioned symmetric Toeplitz systems}, {SIAM J. Sci.
Comput.}, {19-2} (1998), pp.~ 516--529.
%
\bibitem{CN}
{R.~H. Chan and M.~K. Ng}, \textit{Conjugate gradient methods for
Toeplitz systems}, {SIAM Rev.}, {38} (1996), pp.~427--482.
%
\bibitem{mcoseni}
{R.~H. Chan, S. Serra-Capizzano, and C. Tablino-Possio},
\textit{Two-Grid methods for banded linear systems from DCT III
algebra}, {Numer. Linear Algebra Appl.}, {12-2/3} (2005),
pp.~241--249.
%
\bibitem{CSun2}
{Q. Chang, X. Jin, and H. Sun}, \textit{Convergence of the
multigrid method for ill-conditioned block toeplitz systems},
{BIT}, {41-1}(2001), pp.~179--190.
%
\bibitem{mehrmann}
U. Elsner, V. Mehrmann, F. Milde, R.A. R{\"o}mer, and M. Schreiber,
\textit{The Anderson Model of Localization: A Challenge for Modern
Eigenvalue Methods},
{SIAM J. Sci. Comput.}, 20 (1999), pp. 2089--2102.
%
\bibitem{falgout}
R. Falgout, P. Vassilevski, and L. Zikatanov,
\textit{On two-grid convergence estimates}, Numer.
Linear Algebra Appl., 12 (2005), pp. 471--494.
%
\bibitem{FS1}
{G. Fiorentino and S. Serra-Capizzano}, \textit{Multigrid methods
for Toeplitz matrices}, {Calcolo}, {28} (1991),  pp.~283--305.
%
\bibitem{FS2}
{G. Fiorentino and S. Serra-Capizzano}, \textit{Multigrid methods
for symmetric positive definite block Toeplitz matrices
 with nonnegative generating functions},
{SIAM J. Sci. Comput.}, {17-4} (1996), pp.~1068--1081.
%
\bibitem{hack}
{W. Hackbusch}, \textit{Multigrid Methods and Applications},
Springer Verlag, Berlin, 1985.
%
\bibitem{ng-diag-ci}
{M.K. Ho and M.K. Ng}, \textit{Splitting iterations for
circulant-plus-diagonal systems}, Numer. Linear Algebra Appl., {12}
(2005), pp.~779--792.
%
\bibitem{Hu}
{T. Huckle and J. Staudacher}, \textit{Multigrid preconditioning
and Toeplitz matrices}, {Electr. Trans. Numer. Anal.}, {13}
(2002), pp.~81--105.
%
\bibitem{Lundbook}
{J. Lund and K. Bowers}, \textit{Sinc Methods for Quadrature and
Differential Equations}, SIAM Publications, Philadelphia (PA),
1992.
%
\bibitem{super}
{M.~K. Ng, R.~H. Chan, T. Chan, and A. Yip}, \textit{Cosine
Transform Preconditioners for High Resolution Image
Reconstruction}, Linear Algebra Appl., {316} (2000), pp.~89--104.
%
\bibitem{NCT}
{M.~K. Ng,  R.~H. Chan, and W.~C. Tang}, \textit{A fast algorithm
for deblurring models with Neumann boundary conditions}, {SIAM J.
Sci. Comput.}, {21} (1999), pp.~851--866.
%
\bibitem{nst}
{M.~K. Ng, S. Serra-Capizzano, and C. Tablino-Possio},
\textit{Numerical behaviour of multigrid methods for symmetric
Sinc-Galerkin systems}, {Numer. Linear Algebra Appl.}, {12-2/3}
(2005), pp.~261--269.
%
\bibitem{TCS}
 {D. Noutsos, S. Serra-Capizzano, and P. Vassalos},
 \textit{Matrix algebra preconditioners for
 multilevel Toeplitz systems do not insure optimal convergence rate},
 {Theoret. Computer Science}, {315} (2004), pp.~557--579.
%
\bibitem{RStub}
{J. Ruge and K. St\"uben}, \textit{Algebraic Multigrid},
  in Frontiers in Applied Mathematics:
  Multigrid Methods, S. McCormick Ed.,
  {\it SIAM}, Philadelphis (PA), (1987), pp.~73-130.
%
\bibitem{Smulti}
{S. Serra-Capizzano}, \textit{Multi-iterative methods}, {Comput.
Math. Appl.}, {26-4} (1993), pp.~65--87.
%
\bibitem{Sun}
{S. Serra-Capizzano}, \textit{Convergence analysis of two grid
methods for elliptic Toeplitz and PDEs matrix sequences}, {Numer.
Math.}, {92-3} (2002), pp.~433--465.
%
%
\bibitem{model-tau}
{S. Serra-Capizzano}, \textit{A note on anti-reflective boundary
conditions and fast deblurring models}, SIAM J. Sci. Comput.,
{25-3} (2003), pp.~1307--1325.
%
\bibitem{mcirco}
{S. Serra-Capizzano and C. Tablino-Possio}, \textit{Multigrid
methods for multilevel circulant matrices}, {SIAM J. Sci.
Comput.}, {26-1} (2004), pp.~55--85.
%
\bibitem{ST-weighted}
{S. Serra-Capizzano and C. Tablino-Possio}, \textit{A Note on
Algebraic Multigrid Methods for the Discrete Weighted Laplacian},
{Computers $\&$ Mathematics with Applications}, {60-5} (2010),
pp.~1290--1298.
%
%
\bibitem{mcosine-vcycle}
{C. Tablino-Possio}, \textit{V-cycle optimal convergence for
DCT-III matrices},  in {Operator Theory: Advances and
Applications.}, 199 - Numerical Methods for Structured Matrices
and Applications: Georg Heinig memorial volume", D. Bini, V.
Mehrmann, V. Olshevsky, E. Tyrtyshnikov, M. Van Barel Eds.,
Birkhauser Verlag (2010), pp.~377--396.
%
\bibitem{T-LAA-1995} 
{E. Tyrtyshnikov},  \textit{Circulant preconditioners with
unbounded inverses.},  Linear Algebra Appl.,  {216}  (1995),
pp.~1--23.
%
\bibitem{Oost}
{U. Trottenberg, C.~W. Oosterlee, and A. Sch{\"u}ller},
\textit{Multigrid}, Academic Press, London, 2001.
%
\end{thebibliography}
\end{document}